\input amstex

\magnification\magstephalf
\documentstyle{amsppt}

\hsize 5.72 truein
\vsize 7.9 truein
\hoffset .39 truein
\voffset .26 truein
\mathsurround 1.67pt
\parindent 20pt
\normalbaselineskip 13.8truept
\normalbaselines
\binoppenalty 10000
\relpenalty 10000
\csname nologo\endcsname 


\font\bc=cmb10
\font\tenbsy=cmbsy10

\catcode`\@=11

\def\myitem#1.{\item"(#1)."\advance\leftskip10pt\ignorespaces}

\def\qedsymbol{{\mathsurround\z@$\square$}}
\redefine\qed{\relaxnext@\ifmmode\let\next\@qed\else
  {\unskip\nobreak\hfil\penalty50\hskip2em\null\nobreak\hfil
    \qedsymbol\parfillskip\z@\finalhyphendemerits0\par}\fi\next}
\def\@qed#1$${\belowdisplayskip\z@\belowdisplayshortskip\z@
  \postdisplaypenalty\@M\relax#1
  $$\par{\lineskip\z@\baselineskip\z@\vbox to\z@{\vss\noindent\qed}}}
\outer\redefine\beginsection#1#2\par{\par\penalty-250\bigskip\vskip\parskip
  \leftline{\tenbsy x\bf#1. #2}\nobreak\smallskip\noindent}
\outer\redefine\genbeginsect#1\par{\par\penalty-250\bigskip\vskip\parskip
  \leftline{\bf#1}\nobreak\smallskip\noindent}

\def\next{\let\@sptoken= }\def\next@{ }\expandafter\next\next@
\def\@futureletnext#1{\let\nextii@#1\futurelet\next\@flti}
\def\@flti{\ifx\next\@sptoken\let\next@\@fltii\else\let\next@\nextii@\fi\next@}
\expandafter\def\expandafter\@fltii\next@{\futurelet\next\@flti}

\let\zeroindent\z@
\let\savedef@\endproclaim\let\endproclaim\relax 
\define\chkproclaim@{\add@missing\endroster\add@missing\enddefinition
  \add@missing\endproclaim
  \envir@stack\endproclaim
  \edef\endit@{\leftskip\the\leftskip\rightskip\the\rightskip}}
\let\endproclaim\savedef@
\def\thing@{.\enspace\egroup\ignorespaces}
\def\thingi@(#1){ \rm(#1)\thing@}
\def\thingii@\cite#1{ \rm\@pcite{#1}\thing@}
\def\thingiii@{\ifx\next(\let\next\thingi@
  \else\ifx\next\cite\let\next\thingii@\else\let\next\thing@\fi\fi\next}
\def\thing#1#2#3{\chkproclaim@
  \ifvmode \medbreak \else \par\nobreak\smallskip \fi
  \noindent\advance\leftskip#1
  \hskip-#1#3\bgroup\bc#2\unskip\@futureletnext\thingiii@}
\let\savedef@\endproclaim\let\endproclaim\relax 
\def\endit{\endproclaim\endit@\let\endit@\undefined}
\let\endproclaim\savedef@
\def\defn#1{\thing\parindent{Definition #1}\rm}
\def\lemma#1{\thing\parindent{Lemma #1}\sl}
\def\prop#1{\thing\parindent{Proposition #1}\sl}
\def\thm#1{\thing\parindent{Theorem #1}\sl}
\def\cor#1{\thing\parindent{Corollary #1}\sl}

\def\remk#1{\thing\zeroindent{Remark #1}\rm}

\def\narrowthing#1{\chkproclaim@\medbreak\narrower\noindent
  \it\def\next{#1}\def\next@{}\ifx\next\next@\ignorespaces
  \else\bgroup\bc#1\unskip\let\next\narrowthing@\fi\next}
\def\narrowthing@{\@futureletnext\thingiii@}

\def\@cite#1,#2\end@{{\rm([\bf#1\rm],#2)}}
\def\cite#1{\in@,{#1}\ifin@\def\next{\@cite#1\end@}\else
  \relaxnext@{\rm[\bf#1\rm]}\fi\next}
\def\@pcite#1{\in@,{#1}\ifin@\def\next{\@cite#1\end@}\else
  \relaxnext@{\rm([\bf#1\rm])}\fi\next}

\advance\minaw@ 1.2\ex@
\atdef@[#1]{\ampersand@\let\@hook0\let\@twohead0\brack@i#1,\z@,}
\def\brack@{\z@}
\let\@@hook\brack@
\let\@@twohead\brack@
\def\brack@i#1,{\def\next{#1}\ifx\next\brack@
  \let\next\brack@ii
  \else \expandafter\ifx\csname @@#1\endcsname\brack@
    \expandafter\let\csname @#1\endcsname1\let\next\brack@i
    \else \Err@{Unrecognized option in @[}%
  \fi\fi\next}
\def\brack@ii{\futurelet\next\brack@iii}
\def\brack@iii{\ifx\next>\let\next\brack@gtr
  \else\ifx\next<\let\next\brack@less
    \else\relaxnext@\Err@{Only < or > may be used here}
  \fi\fi\next}
\def\brack@gtr>#1>#2>{\setboxz@h{$\m@th\ssize\;{#1}\;\;$}%
 \setbox@ne\hbox{$\m@th\ssize\;{#2}\;\;$}\setbox\tw@\hbox{$\m@th#2$}%
 \ifCD@\global\bigaw@\minCDaw@\else\global\bigaw@\minaw@\fi
 \ifdim\wdz@>\bigaw@\global\bigaw@\wdz@\fi
 \ifdim\wd@ne>\bigaw@\global\bigaw@\wd@ne\fi
 \ifCD@\enskip\fi
 \mathrel{\mathop{\hbox to\bigaw@{$\ifx\@hook1\lhook\mathrel{\mkern-9mu}\fi
  \setboxz@h{$\displaystyle-\m@th$}\ht\z@\z@
  \displaystyle\m@th\copy\z@\mkern-6mu\cleaders
  \hbox{$\displaystyle\mkern-2mu\box\z@\mkern-2mu$}\hfill
  \mkern-6mu\mathord\ifx\@twohead1\twoheadrightarrow\else\rightarrow\fi$}}%
 \ifdim\wd\tw@>\z@\limits^{#1}_{#2}\else\limits^{#1}\fi}%
 \ifCD@\enskip\fi\ampersand@}
\def\brack@less<#1<#2<{\setboxz@h{$\m@th\ssize\;\;{#1}\;$}%
 \setbox@ne\hbox{$\m@th\ssize\;\;{#2}\;$}\setbox\tw@\hbox{$\m@th#2$}%
 \ifCD@\global\bigaw@\minCDaw@\else\global\bigaw@\minaw@\fi
 \ifdim\wdz@>\bigaw@\global\bigaw@\wdz@\fi
 \ifdim\wd@ne>\bigaw@\global\bigaw@\wd@ne\fi
 \ifCD@\enskip\fi
 \mathrel{\mathop{\hbox to\bigaw@{$%
  \setboxz@h{$\displaystyle-\m@th$}\ht\z@\z@
  \displaystyle\m@th\mathord\ifx\@twohead1\twoheadleftarrow\else\leftarrow\fi
  \mkern-6mu\cleaders
  \hbox{$\displaystyle\mkern-2mu\copy\z@\mkern-2mu$}\hfill
  \mkern-6mu\box\z@\ifx\@hook1\mkern-9mu\rhook\fi$}}%
 \ifdim\wd\tw@>\z@\limits^{#1}_{#2}\else\limits^{#1}\fi}%
 \ifCD@\enskip\fi\ampersand@}


\define\today{\number\day\ \ifcase\month\or
  January\or February\or March\or April\or May\or June\or
  July\or August\or September\or October\or November\or December\fi
  \ \number\year}
\def\pr@m@s{\ifx'\next\let\nxt\pr@@@s \else\ifx^\next\let\nxt\pr@@@t
  \else\let\nxt\egroup\fi\fi \nxt}

\define\widebar#1{\mathchoice
  {\setbox0\hbox{\mathsurround\z@$\displaystyle{#1}$}\dimen@.1\wd\z@
    \ifdim\wd\z@<.4em\relax \dimen@ -.16em\advance\dimen@.5\wd\z@ \fi
    \ifdim\wd\z@>2.5em\relax \dimen@.25em\relax \fi
    \kern\dimen@ \overline{\kern-\dimen@ \box0\kern-\dimen@}\kern\dimen@}%
  {\setbox0\hbox{\mathsurround\z@$\textstyle{#1}$}\dimen@.1\wd\z@
    \ifdim\wd\z@<.4em\relax \dimen@ -.16em\advance\dimen@.5\wd\z@ \fi
    \ifdim\wd\z@>2.5em\relax \dimen@.25em\relax \fi
    \kern\dimen@ \overline{\kern-\dimen@ \box0\kern-\dimen@}\kern\dimen@}%
  {\setbox0\hbox{\mathsurround\z@$\scriptstyle{#1}$}\dimen@.1\wd\z@
    \ifdim\wd\z@<.28em\relax \dimen@ -.112em\advance\dimen@.5\wd\z@ \fi
    \ifdim\wd\z@>1.75em\relax \dimen@.175em\relax \fi
    \kern\dimen@ \overline{\kern-\dimen@ \box0\kern-\dimen@}\kern\dimen@}%
  {\setbox0\hbox{\mathsurround\z@$\scriptscriptstyle{#1}$}\dimen@.1\wd\z@
    \ifdim\wd\z@<.2em\relax \dimen@ -.08em\advance\dimen@.5\wd\z@ \fi
    \ifdim\wd\z@>1.25em\relax \dimen@.125em\relax \fi
    \kern\dimen@ \overline{\kern-\dimen@ \box0\kern-\dimen@}\kern\dimen@}%
  }

\catcode`\@\active

\let\PVstyle=d 

\loadeufm

\font\tenscr=rsfs10 
\font\sevenscr=rsfs7 
\font\fivescr=rsfs5 
\skewchar\tenscr='177 \skewchar\sevenscr='177 \skewchar\fivescr='177
\newfam\scrfam \textfont\scrfam=\tenscr \scriptfont\scrfam=\sevenscr
\scriptscriptfont\scrfam=\fivescr
\define\scr#1{{\fam\scrfam#1}}
\let\Cal\scr

\outer\redefine\subhead#1\par{\par\medbreak\leftline{\bf #1}\smallskip\noindent}

\loadbold

\let\0\relax 
\mathchardef\idot="202E

\define\bir{{\text{bir}}}
\define\codim{\operatorname{codim}}
\define\exc{{\text{exc}}}
\define\lvol{\operatorname{vol}}
\define\Nev{\operatorname{Nev}}
\define\Nevbir{\Nev_{\operatorname{bir}}}
\define\Pic{\operatorname{Pic}}

\define\Spec{\operatorname{Spec}}
\define\sq#1{\ifx#1([\else\ifx#1)]\else\message{invalid use of "sq"}\fi\fi}
\define\Supp{\operatorname{Supp}}

\define\restrictedto#1{\big|_{#1}}

\catcode`\@=11
\define\bibdef#1#2#3{\expandafter\def\csname by.#1\endcsname{#2}
  \expandafter\def\csname ba.#1\endcsname{#3}}
\define\bib@chkdoyear#1{\expandafter\ifx\csname by.#1\endcsname\relax
  \message{*** Warning: Undefined bibliography tag #1 ***}\fi
  \csname by.#1\endcsname}

\define\@citepguts#1{[\bf\csname ba.#1\endcsname\ \bib@chkdoyear{#1}\rm]}
\define\@citep#1,#2\end@{{\rm(\@citepguts{#1},#2)}}
\define\citep#1{\in@,{#1}\ifin@\def\next{\@citep#1\end@}\else
  \relaxnext@{\rm\@citepguts{#1}}\fi\next}

\redefine\@pcite#1{\in@,{#1}\ifin@\def\next{\@citep#1\end@}\else
  \relaxnext@{\rm(\@citepguts{#1})}\fi\next}

\define\@citetauth#1{\csname ba.#1\endcsname}
\define\@citetyear#1{[\bf\bib@chkdoyear{#1}\rm]}
\define\@citet#1,#2\end@{{\rm\@citetauth{#1} (\@citetyear{#1},#2)}}
\define\citet#1{\in@,{#1}\ifin@\def\next{\@citet#1\end@}\else
  \relaxnext@{\rm\@citetauth{#1} \@citetyear{#1}}\fi\next}

\define\auth#1{\par\bigbreak\leftline{\bf #1}\par}

\define\refer#1{\par
  \hangafter 1%
  \hangindent \bibindent
  \noindent \hbox to \bibindent{\hskip\bibitemindent \bib@chkdoyear{#1}\hfil}%
  \ignorespaces}

\catcode`\@\active

\bibdef{a}{2011}{Autissier}
\bibdef{ha77}{1977}{Hartshorne}
\bibdef{hl17}{2017pre}{Heier and Levin}
\bibdef{JDB}{2011}{Jothilingham, et al\.}
\bibdef{la_fdg}{1983}{Lang}
\bibdef{laz04}{2004}{Lazarsfeld}
\bibdef{rv1}{2020}{Ru and Vojta}
\bibdef{rv2}{2020pre}{Ru and Vojta}
\bibdef{rw17}{2017}{Ru and Wang}
\bibdef{rw20p}{2020pre}{Ru and Wang}
\bibdef{sil87}{1987}{Silverman}
\bibdef{stks}{2020}{Stacks project authors}
\bibdef{voj89}{1989}{Vojta}
\bibdef{voj11}{2011}{Vojta}
\bibdef{y04}{2004}{Yamanoi}

\topmatter
\title Birational Nevanlinna Constants, Beta Constants, and
  Diophantine Approximation to Closed Subschemes\endtitle
\rightheadtext{Birational Nevanlinna constants and closed subschemes}
\author Paul Vojta\endauthor
\affil University of California, Berkeley\endaffil
\address Department of Mathematics, University of California,
  970 Evans Hall\quad\#3840, Berkeley, CA \ 94720-3840\endaddress
\date 19 September 2020 \enddate

\abstract
In an earlier paper (joint with Min Ru), we proved a result on
diophantine approximation to Cartier divisors, extending a 2011 result of
P. Autissier.  This was recently extended to certain closed subschemes
(in place of divisors) by Ru and Wang.
In this paper we extend this result to a broader class of closed subschemes.
We also show that some notions of $\beta(\Cal L,D)$ coincide,
and that they can all be evaluated as limits.
\endabstract
\endtopmatter

\document

Let $k$ be either a number field or the field $\Bbb C$ of complex numbers,
and let $X$ be a complete variety over $k$ (see Section \01 for detailed
definitions).
We recall the following from \citep{rv1, Def.~1.9 and ``General Theorem''}.

\defn{\00.1}  Let $\Cal L$ be a big line sheaf on $X$ and let $D$ be a
nonzero effective Cartier divisor on $X$.  Then
$$\beta(\Cal L, D)
  = \liminf_{N\to\infty} \frac{\sum_{m=1}^\infty h^0(X,\Cal L^N(-mD))}
    {Nh^0(X,\Cal L^N)}\;.\tag\00.1.1$$
(In this paper $\Cal L^N$ always means $\Cal L^{\otimes N}$, the tensor
power of $N$ copies of $\Cal L$.)
\endit

\thm{\00.2}  Let $k$ and $X$ be as above, let $\Cal L$ be a big line sheaf
on $X$, and let $D_1,\dots,D_q$ be nonzero effective Cartier divisors on $X$
that intersect properly (i.e., for any nonempty $I\subseteq\{1,\dots,q\}$
and any $x\in\bigcap_{i\in I}\Supp D_i$, the divisors $D_i$, $i\in I$
are locally generated near $x$ by a regular sequence in $\Cal O_{X,x}$).
\roster
\myitem a. {\bc (Arithmetic part)}  Assume that $k$ is a number field,
and let $S$ be a finite set of places of $k$.  Then, for all $\epsilon>0$,
there is a proper Zariski-closed subset $Z$ of $X$ such that the inequality
$$\sum_{i=1}^q \beta(\Cal L,D_i) m_S(D_i,x)
    \le (1+\epsilon) h_{\Cal L}(x) + O(1)\tag\00.2.1$$
holds for all points $x\in X(k)\setminus Z$.
\myitem b. {\bc (Analytic part)}  Assume that $k=\Bbb C$.  Then,
for all $\epsilon>0$, there is a proper Zariski-closed subset $Z$ of $X$
such that the inequality
$$\sum_{i=1}^q \beta(\Cal L,D_i) m_f(D_i,r)
    \le_\exc (1+\epsilon) T_{f, \Cal L}(r)\tag\00.2.2$$
holds for all holomorphic mappings $f\:\Bbb C\to X$ whose image is
not contained in $Z$.  The subscript ``exc'' means that the inequality holds
for all $r\in(0,\infty)$ outside of a set of finite Lebesgue measure.
\endroster
\endit

\remk{\00.3}  Part (b) in this theorem has been strengthened so that it
applies to all $f$ whose image is not contained in $Z$, whereas
in \citep{rv1} $f$ was required to have Zariski-dense image.
This version can be obtained by revising the statements of Thm.~2.7
(see Remark 2.8), Thm.~2.11, Thm.~1.4, and the Main Theorem
of \citep{rv1} accordingly, where $Z$ depends on $\epsilon$ only in
the last two theorems.
\endit

The purpose of this paper is to generalize Theorem \00.2 to replace the
divisors $D_i$ with proper closed subschemes $Y_i$.

Upon circulating an early version of this paper, I was informed that
\citet{rw20p} had already proved a version of Theorem \00.2 for
closed subschemes.  However, the version presented here is somewhat
more general.

For both the work of Ru and Wang and the present paper, extending
Theorem \00.2 to closed subschemes involves defining what it means for
the subschemes $Y_i$ to intersect properly.  In both cases this is
done using regular sequences---see Remark \02.11 and Definition \03.1.
However, the details of this definition are different in the two papers,
and this is the main difference between them.

For example, if $X$ is Cohen--Macaulay (e.g., if it is nonsingular),
then the $Y_i$ intersect properly, in the sense of the present paper,
if and only if (i) at each intersection point $x$, each of the $Y_i$
passing through $x$ is generated by monomials in the elements of some regular
sequence in the local ring, and (ii) the subschemes $Y_i$ are in general
position (in other words, they intersect properly in the sense of
intersection theory).
See Definitions \02.4 and \02.9.  This condition is only needed at points
where two or more of the $Y_i$ intersect, leading to a definition that
they ``weakly intersect properly'' (Definition \03.1c).
The definition of Ru and Wang uses the stronger condition
that the ideals are generated by the actual elements of a regular sequence.
In particular, (in the Cohen--Macaulay case) their result requires each $Y_i$
to be a local complete intersection as a scheme, but this paper relaxes
this condition somewhat---see Remark \02.11.

The generalization of Theorem \00.2 to be proved here is stated below
as Theorem \00.9.  This statement also describes the main theorem of
\citet{rw20p}, except that it is relative to the notion of proper intersection
described in Remark \02.11 instead of Definition \03.1.

\citet{hl17} have also proved a diophantine theorem on approximation to
proper closed subschemes.  In their theorem, closed subschemes
of codimension $r$ may be repeated up to $r$ times.  In this paper, as well
as in the paper of Ru and Wang, however, subschemes may not be repeated.
Instead, $\beta(\Cal L,Y_i)$ is usually larger for such subschemes,
as is the case for linear subspaces of projective space (Proposition \09.2).

The definition of $\beta(\Cal L,Y)$ for a proper closed subscheme $Y$ of $X$
is a straightforward extension of (\00.1.1):

\defn{\00.4} \cite{rw20p, Def.~1.2}  Let $\Cal L$ be a big line sheaf on $X$,
let $Y$ be a nonempty proper closed subscheme of $X$, and let $\Cal I$ be
the sheaf of ideals corresponding to $Y$.  Then
$$\beta(\Cal L, Y)
  = \liminf_{N\to\infty} \frac{\sum_{m=1}^\infty h^0(X,\Cal L^N\otimes\Cal I^m)}
    {Nh^0(X,\Cal L^N)}\;.\tag\00.4.1$$
\endit

\remk{\00.5}  A closely related definition was given by \citet{rw17, Def.~1.1}:
$$\beta_{\Cal L,Y}
  = \liminf_{N\to\infty} \frac{\sum_{m=1}^\infty h^0(W,\pi^{*}\Cal L^N(-mE))}
    {Nh^0(X,\Cal L^N)}\;,\tag\00.5.1$$
where $\pi\:W\to X$ is the blowing-up of $X$ along $Y$ and $E$ is the
exceptional divisor, so in particular the two definitions coincide when
$Y$ is an effective Cartier divisor.  In fact, they coincide for all $Y$;
see \citep{rw20p, Remark~1.3} when $X$ is Cohen--Macaulay and $Y$ is a
local complete intersection, or Corollary \05.9 for the general case.
\endit

Another major goal of this paper is to show that these three beta constants
all coincide (Corollary \05.9).

\remk{\00.6}  There is a ``birational'' version of Definition \00.1,
in which $D$ is replaced by a Cartier b-divisor $\bold D$.
This constant is denoted $\beta(\Cal L,\bold D)$; see Definition \05.5.
Since a proper closed subscheme is a special case of a b-divisor
(see Definition \05.7), this leads to a constant $\beta(\Cal L,\bold Y)$
defined by a slightly different limit.  This, as it turns out, has the
same value as $\beta(\Cal L,Y)$ and $\beta_{\Cal L,Y}$---see Corollary \05.9.
Note also that in Definitions \00.1, \00.5, \05.5, and \00.4,
the $\liminf$ can be replaced by a limit whenever $\operatorname{char} k=0$.
This is proved in Section \010.
\endit

\remk{\00.7}  The proof of Theorem \00.9 uses $\beta(\Cal L,Y)$.
This is because the Autissier property (see below) is not preserved
by blowing up, so we work on $Y$.
\endit

\remk{\00.8}  It is possible to let $\bold D$ be an $\Bbb R$\snug-Cartier
b-divisor in the definition of $\beta(\Cal L,\bold D)$.  We have not done that
here, though, as it would not provide any benefit (so far), but would
involve additional complexity.
\endit

One defines Weil functions relative to proper closed subschemes $Y$ on $X$
by blowing up $X$ along $Y$ to obtain a Cartier divisor on the blow-up;
see for example \citep{sil87, 2.2} or \citep{y04, 2.2}, in combination with
\citep{sil87, Thm.~2.1(h)}.  These can then be used to define proximity
and counting functions for $Y$.  For details see Section \08.

With these definitions, the main theorem of this paper is as follows.

\thm{\00.9}  Let $X$ be a complete variety over a field $k$, let $\Cal L$
be a big line sheaf on $X$, and let $Y_1,\dots,Y_q$ be
proper closed subschemes of $X$ that weakly intersect properly
(see Definition \03.1).
\roster
\myitem a. {\bc (Arithmetic part)}  Assume that $k$ is a number field,
and let $S$ be a finite set of places of $k$.  Then, for all $\epsilon>0$
and all $C\in\Bbb R$, there is a proper Zariski-closed subset $Z$ of $X$
such that the inequality
$$\sum_{i=1}^q \beta(\Cal L,Y_i) m_S(Y_i,x)
    \le (1+\epsilon) h_{\Cal L}(x) + C\tag\00.9.1$$
holds for all points $x\in X(k)\setminus Z$.
\myitem b. {\bc (Analytic part)}  Assume that $k=\Bbb C$.  Then,
for all $\epsilon>0$, there is a proper Zariski-closed subset $Z$ of $X$
such that the inequality
$$\sum_{i=1}^q \beta(\Cal L,Y_i) m_f(Y_i,r)
    \le_\exc (1+\epsilon) T_{f,\Cal L}(r)\tag\00.9.2$$
holds for all holomorphic mappings $f\:\Bbb C\to X$ whose image is
not contained in $Z$.
\endroster
\endit

Diophantine inequalities for closed subschemes have already been obtained
by other authors.  For example, \citet{rw17} proved the inequality
$$\sum_{i=1}^q m_S(Y_i,x)
  \le \left(\ell\max_i\beta(\Cal L,Y_i)^{-1}+\epsilon\right) h_{\Cal L}(x)\;,$$
where at most $\ell$ of the $Y_i$ have nonempty intersection.
This overlaps with our results here, but is not fully implied by our
Theorem \00.9 because the latter theorem requires that the $Y_i$
weakly intersect properly, but Ru and Wang only need the condition
involving $\ell$.

\citet{hl17} also have an inequality involving closed subschemes.
Their theorem again has weaker conditions on the $Y_i$ (and in fact it
allows some $Y_i$ of codimension $>1$ to be repeated).  It is harder to
compare their theorem to ours since their theorem involves Seshadri
constants.  Not much is known about how Seshadri constants compare
with $\beta(\Cal L,Y_i)$.

Theorem \00.9 will be proved by splitting it up into two theorems,
involving a property due originally to \citet{a, Lemme 3.3};
see also Lemma \02.3.  This will be expressed by saying that closed
subschemes $Y_1,\dots,Y_q$ {\bc have the Autissier property};
see Definition \03.2.

These two theorems are the following.

\thm{\00.10}  Let $X$ be a complete variety over a field $k$, and
let $Y_1,\dots,Y_q$ be proper closed subschemes of $X$.  If $Y_1,\dots,Y_q$
weakly intersect properly, then they have the Autissier property.
\endit

\thm{\00.11}  Let $X$ be a complete variety over a field $k$
of characteristic $0$, let $\Cal L$ be a big line sheaf on $X$,
and let $Y_1,\dots,Y_q$ be proper closed subschemes of $X$
that have the Autissier property.
Then (depending on $k$) part (a) or (b) of Theorem \00.9 holds.
\endit


It is clear that the conjunction of these theorems implies Theorem \00.9.

The outline of the paper is as follows.  Section \01 briefly gives some
fundamental definitions used in the paper.  Sections \02 and \03 give a
version of Autissier's lemmas on ideals associated to saturated subsets
of $\Bbb N^r$, in the local and global cases, respectively, leading up to
the proof of Theorem \00.10 (see Proposition \03.3) and
the proof in Section \04 that his function $N(\bold t)$ is convex.
These constitute the key insight of the paper.  Section \05 develops
the machinery that will be used to work with the ideal sheaves associated
to the closed subschemes $Y_i$ in the paper.  Section \06 adapts work
of \citet{a}, as modified by \citet{rv1}, to the current context,
finishing the technical parts of the proof.  Section \07 introduces a
birational Nevanlinna constant for $\Bbb R$\snug-Cartier b-divisors,
which is then used in Section \08 to finish the proof of Theorem \00.9.
Section \09 explores the special case of linear subvarieties of $\Bbb P^n$.
Finally, Section \010 gives a detailed proof of the fact that the limits infima
in (\00.1.1), (\00.4.1), (\00.5.1), and (\05.5.1) (the definitions of
$\beta(\Cal L,D)$, $\beta(\Cal L,Y)$, $\beta_{\Cal L,Y}$,
and $\beta(\Cal L,\bold D)$, respectively) can be replaced by limits
(in characteristic $0$).

I thank Min Ru for suggesting the idea of extending the Main Theorem
of \citep{rv1} to subschemes.

\beginsection{\01}{Basic Notation and Conventions}

The basic notations of this paper follow those of \citep{rv1} and \citep{rv2}.

In this paper $\Bbb N=\{0,1,2,\dots\}$.  Also $\Bbb Z_{>0}=\{1,2,3,\dots\}$,
$\Bbb R_{\ge0}=\{x\in\Bbb R:x\ge0\}$, etc.

A {\bc variety} over a field $k$ is an integral scheme, separated and
of finite type over $k$.  A {\bc morphism} of varieties over $k$ is
a morphism of schemes over $k$.

Subschemes will always be assumed to be closed and proper (i.e.,
not the whole scheme).

\beginsection{\02}{A Property of Autissier}

This section extends \citep{a, Lemme 3.3} to accommodate subschemes of
higher codimension.

This lemma motivates a definition of a property of subschemes,
which basically says that they satisfy the conclusion of this lemma.
This property will be called the {\it Autissier property\/};
see Definitions \02.12 and \03.2.  The entire remainder of the proof of
Theorem \00.9 hinges on this property.

Throughout this section, $A$ is a noetherian local ring.

We start by recalling some definitions and a lemma of \citet{a}.

\defn{\02.1}  Let $r\in\Bbb Z_{>0}$.  A subset $N$ of $\Bbb N^r$ is
{\bc saturated} if it is nonempty and if $N\supseteq \bold a+\Bbb N^r$
for all $\bold a\in N$.
\endit

\defn{\02.2}  Let $\phi_1,\dots,\phi_r\in A$ with $r>0$, and let $N$ be a
saturated subset of $\Bbb N^r$.  Then $\Cal I(N)$ is the ideal of $A$
generated by the set $\{\phi_1^{b_1}\dots\phi_r^{b_r}:\bold b\in N\}$.
\endit

The key fact about this definition is the following lemma due to Autissier.

\lemma{\02.3} \cite{a, Lemme 3.3}  Let $\phi_1,\dots,\phi_r$ ($r>0$) be
a regular sequence in $A$, and let $N_1$ and $N_2$ be saturated subsets
of $\Bbb N^r$.  Then
$$\Cal I(N_1\cap N_2) = \Cal I(N_1) \cap \Cal I(N_2)\;.$$
\endit

Now we carry the above over to the situation of ideals in $A$.

\defn{\02.4}  Let $I$ be an ideal of $A$ and let $\phi_1,\dots,\phi_r$ be
a sequence of elements of $A$.  Then $I$ is of {\bc monomial type}
with respect to $\phi_1,\dots,\phi_r$ if $r>0$ and $I=\Cal I(N)$
(taken relative to $\phi_1,\dots,\phi_r$)
for some saturated subset $N$ of $\Bbb N^r$.
\endit

Note that if $I$ is of monomial type with respect to some sequence
$\phi_1,\dots,\phi_r$, then so is $I^n$ for all $n\in\Bbb N$.
This is immediate from the following lemma.

\lemma{\02.5}  Let $r\in\Bbb Z_{>0}$ and let $N$ be a saturated subset
of $\Bbb N^r$.  For all $n\in\Bbb N$ let
$$nN = \cases \Bbb N^r & \text{if $n=0$;} \\
  \{\bold b_1+\dots+\bold b_n:\bold b_1,\dots,\bold b_n\in N\}
    & \text{if $n>0$}\;.
  \endcases\tag\02.5.1$$
(When $n>0$ this is the Minkowski sum of $N$ with itself $n$ times.)  Then
\roster
\myitem a.  $nN$ is saturated for all $n$;
\myitem b.  $\Cal I(N)^n=\Cal I(nN)$ for all $n$; and
\myitem c.  $nN\subseteq mN$ for all $n\ge m\ge 0$.
\endroster
\endit

\demo{Proof}  Left to the reader.\qed\enddemo

As a counterpart to Definition \02.2, but with closed subschemes in place of
Cartier divisors, we have the following.

\defn{\02.6}  Let $q\in\Bbb Z_{>0}$, let $I_1,\dots,I_q$ be ideals in $A$,
and let $N$ be a saturated subset of $\Bbb N^q$.  Then $\Cal J(N)$ is
the ideal of $A$ defined by
$$\Cal J(N) = \sum_{\bold b\in N} I_1^{b_1}\dotsm I_q^{b_q}\;.$$
\endit

This can be expressed in terms of $\Cal I(\cdot)$ as follows.

\defn{\02.7}  Let $q\in\Bbb Z_{>0}$.  For each $i=1,\dots,q$
let $M_i$ be a saturated subset of $\Bbb N^{r_i}$ with $r_i\in\Bbb Z_{>0}$.
For all saturated subsets $N$ of $\Bbb N^q$, we then define
$$M(N) = \bigcup_{\bold c\in N} c_1M_1\times\dots\times c_qM_q\;.\tag\02.7.1$$
This is a saturated subset of $\Bbb N^r$, where $r=r_1+\dots+r_q$.
\endit

\lemma{\02.8}  Let $q\in\Bbb Z_{>0}$.  For each $i=1,\dots,q$,
let $M_i$ be a saturated subset of $\Bbb N^{r_i}$ and let $I_i=\Cal I(M_i)$,
taken relative to a nonempty sequence $\phi_{i1},\dots,\phi_{ir_i}$ in $A$.
Let $N$ be a saturated subset of $\Bbb N^q$.  Then
$$\Cal J(N) = \Cal I(M(N))\;,$$
where $\Cal J(N)$ is taken with respect to $I_1,\dots,I_q$
and $\Cal I(M(N))$ is taken with respect to the sequence
$$\phi_{11},\dots,\phi_{1r_1},\dots,\phi_{q1},\dots,\phi_{qr_q}\;.\tag\02.8.1$$
\endit

\demo{Proof}  This is immediate from Definitions \02.2, \02.6, and \02.7.
See also \citep{rw20p, Lemma~3.3}.\qed
\enddemo

We can now state the main definitions and main result of this section.

\defn{\02.9}  Let $I_1,\dots,I_q$ be ideals of $A$, with $q\in\Bbb N$.
Then $I_1,\dots,I_q$ {\bc intersect properly} if (i) for each $i=1,\dots,q$
there is a nonempty regular sequence $\phi_{i1},\dots,\phi_{ir_i}$ in $A$ such
that $I_i$ is of monomial type with respect to $\phi_{i1},\dots,\phi_{ir_i}$;
and (ii) the combined sequence (\02.8.1) is a regular sequence.
\endit

\remk{\02.10}  Since the length of the sequence (\02.8.1) is at most $\dim A$,
we must have $q\le\dim A$ whenever $I_1,\dots,I_q$ intersect properly.
\endit

\remk{\02.11}  \citet{rw20p} say that $I_1,\dots,I_q$ intersect properly if,
in the notation of Definition \02.9, $I_i$ is generated by
$\phi_{i1},\dots,\phi_{ir_i}$ for all $i$ (and (\02.8.1) is a regular sequence).
In other words, this is the special case of Definition \02.9 in which
the subset $N$ of Definition \02.4 equals $\Bbb N^r\setminus\{\bold 0\}$.
One then obtains their notion of subschemes $Y_1,\dots,Y_q$
intersecting properly by using this definition in place of Definition \02.9
in Definition \03.1.

As an example, assume that $X$ contains an open subset isomorphic to
$\Bbb A^2_k=\allowmathbreak\Spec k[x,y]$.  Then $x,y$ is a regular sequence
in the local
ring at $(0,0)$, so $I=(x,y)$ satisfies the hypotheses of Ru and Wang's
theorem (and it is also of monomial type with respect to $x$ and $y$).
The ideal $(x^3,xy,y^2)$, however, does not satisfy their condition,
but it is of monomial type in $x,y$, so it can be handled by Theorem \00.9.
\endit

\defn{\02.12}  Let $I_1,\dots,I_q$ be ideals in $A$.
We say that they {\bc have the Autissier property} if
$$\Cal J(N\cap N') = \Cal J(N) \cap \Cal J(N')\tag\02.12.1$$
for all saturated subsets $N$ and $N'$ of $\Bbb N^q$.
\endit

\prop{\02.13}  Let $I_1,\dots,I_q$ be ideals in $A$.  If they
intersect properly, then they have the Autissier property.
\endit

\demo{Proof}  By Lemmas \02.8 and \02.3, we immediately reduce to showing that
$$M(N\cap N') = M(N) \cap M(N')\;.\tag\02.13.1$$

To prove this, we first need some basic facts on the product ordering
on $\Bbb N^q$.

Recall that the product ordering on $\Bbb N^q$ is defined by $\bold a\le\bold b$
if and only if $a_i\le b_i$ for all $i=1,\dots,q$.  This ordering is a lattice;
in particular, for any $\bold a,\bold b\in\Bbb N^q$, the join, or
least upper bound, of $\bold a$ and $\bold b$ is the element
$\bold a\vee\bold b=\bold c\in\Bbb N^q$ defined by $c_i=\max\{a_i,b_i\}$
for all $i$.

Now we note that if $N$ and $N'$ are saturated subsets of $\Bbb N^q$, then
$$\{\bold c\vee\bold c' : \bold c\in N, \bold c'\in N'\} = N\cap N'\;.
  \tag\02.13.2$$
Indeed, the inclusion ``$\supseteq$'' is immediate by taking $\bold c'=\bold c$
for all $\bold c\in N\cap N'$.  Conversely, if $\bold c''=\bold c\vee\bold c'$
with $\bold c\in N$ and $\bold c'\in N'$, then $\bold c''\in N$
and $\bold c''\in N'$ because $N$ and $N'$ are saturated (respectively);
hence $\bold c''\in N\cap N'$.

Then, by (\02.7.1), distributivity of $\cap$ over $\cup$, compatibility of
intersection and product, Lemma \02.5c, (\02.13.2), and (\02.7.1) again,
we have
$$\split M(N) \cap M(N')
  &= \left(\bigcup_{\bold c\in N} c_1M_1\times\dots\times c_qM_q\right)
    \cap
    \left(\bigcup_{\bold c'\in N'} c_1'M_1\times\dots\times c_q'M_q\right) \\
  &= \bigcup\Sb \bold c\in N\\\bold c'\in N'\endSb
    \bigl((c_1M_1\times\dots\times c_qM_q)
      \cap(c_1'M_1\times\dots\times c_q'M_q)\bigr) \\
  &= \bigcup_{\bold c,\bold c'} \bigl((c_1M_1\cap c_1'M_1)
    \times\dots\times(c_qM_q\cap c_q'M_q)\bigr) \\
  &= \bigcup_{\bold c,\bold c'} \max\{c_1,c_1'\}M_1
    \times\dots\times\max\{c_q,c_q'\}M_q \\
  &= \bigcup_{\bold c''\in N\cap N'} c_1''M_1\times\dots\times c_q''M_q \\
  &= M(N\cap N')\;.\endsplit$$
This gives (\02.13.1).\qed
\enddemo

Turning to consequences of the Autissier property, in the local setting
we only need the following.

\prop{\02.14} (\citep{a, Remarque 3.4} and \citep{rv1, Remark 6.3})
Let $q\in\Bbb Z_{>0}$, let
$$\square = \Bbb R_{\ge0}^q\setminus\{\bold 0\}\;,\tag\02.14.1$$
and for all $\bold t\in\square$ and all $x\in\Bbb R_{\ge0}$ let
$$N(\bold t,x)
  = \{\bold b\in\Bbb N^q:t_1b_1+\dots+t_qb_q\ge x\}\;.\tag\02.14.2$$
Let $I_1,\dots,I_q$ be ideals in $A$ that have the Autissier property.  Then
$$\Cal J(N(\bold t,x)) \cap \Cal J(N(\bold u,y))
  \subseteq \Cal J(N(\lambda\bold t+(1-\lambda)\bold u,\lambda x+(1-\lambda)y))
  \tag\02.14.3$$
for all $\bold t,\bold u\in\square$, all $x,y\in\Bbb R_{\ge0}$, and
all $\lambda\in[0,1]$.
\endit

\demo{Proof}  This is immediate from Definition \02.12 and the observation that
$$N(\bold t,x) \cap N(\bold u,y)
  \subseteq N(\lambda\bold t+(1-\lambda)\bold u,\lambda x+(1-\lambda)y)\;.\qed$$
\enddemo

\remk{\02.15}  An interesting theory of regular sequences of ideals
has been developed by \citet{JDB}.  In this theory, ideals $I_1,\dots,I_q$
of $A$ are said to be a {\bc regular sequence of ideals} if all of them are
nonzero and proper, and if
$$(I_1+\dots+I_j)\cap I_{j+1} = (I_1+\dots+I_j) \idot I_{j+1}$$
for all $j=1,\dots,q-1$.  This extends the definition of a regular sequence
of elements of a local ring, in the sense that a sequence $(x_1),\dots,(x_q)$
of principal ideals in $A$ is regular if and only if the elements
$x_1,\dots,x_q$ form a regular sequence.

Although it was very tempting to write this paper using the concept of
regular sequences of ideals, ultimately we decided not to.  This was because
many of the results of \citep{JDB} (e.g., Theorem 1) assumed that $A$
was a regular local ring; in addition, there were other difficulties in trying
to rewrite the proof of \citep{a, Lemme~6.2} directly in terms of a
regular sequence of ideals.
\endit

\beginsection{\03}{The Autissier Property of Subschemes}

This brief section carries over Definitions \02.9 and \02.12 and
Proposition \02.13 to the case of subschemes.

First we start with the definitions.

Throughout this section, $X$ is a complete variety over a field $k$ and
$Y_1,\dots,Y_q$ are proper closed subschemes of $X$.

\defn{\03.1}  Let $\Cal I_1,\dots,\Cal I_q$ be the ideal sheaves
that correspond to $Y_1,\dots,Y_q$, respectively.
\roster
\myitem a.  We say that $Y_1,\dots,Y_q$ {\bc intersect properly} at
a point $P\in X$ if the subsequence of proper ideals in
the sequence $(\Cal I_1)_P,\dots,(\Cal I_q)_P$ of ideals of the local ring
$\Cal O_{X,P}$ intersect properly (in the sense of Definition \02.9).
(If the subsequence is trivial, i.e., if $P\notin\bigcup Y_i$, then
this is vacuously true.)
\myitem b.  We say that $Y_1,\dots,Y_q$ {\bc intersect properly}
if $Y_1,\dots,Y_q$ intersect properly at all points of $X$.
\myitem c.  We say that $Y_1,\dots,Y_q$ {\bc weakly intersect properly}
if they intersect properly at all $P\in\bigcup_{i\ne j} (Y_i\cap Y_j)$.
\endroster
\endit

Clearly, if $Y_1,\dots,Y_q$ intersect properly, then they also
weakly intersect properly.

\defn{\03.2}  Let $\Cal I_1,\dots,\Cal I_q$ be as in Definition \03.1.
\roster
\myitem a.  Let $P\in X$, and let $j_1,\dots,j_r$ be the subsequence of
$1,\dots,q$ consisting of those $j$ such that $P\in Y_j$.  We say that
$Y_1,\dots,Y_q$ {\bc have the Autissier property} at $P$ if
$$\Cal J(N\cap N') = \Cal J(N) \cap \Cal J(N')\tag\03.2.1$$
for all saturated subsets $N$ and $N'$ of $\Bbb N^r$, where $\Cal J$
is taken with respect to the sequence $(\Cal I_{j_1})_P,\dots,(\Cal I_{j_r})_P$
of (proper) ideals of $\Cal O_{X,P}$.  (This is equivalent to saying that
$(\Cal I_{j_1})_P,\dots,(\Cal I_{j_r})_P$ have the Autissier property
as in Definition \02.12).
\myitem b.  We say that $Y_1,\dots,Y_q$ {\bc have the Autissier property}
if they have the Autissier property at all $P\in X$.
\endroster
\endit

Corresponding to Proposition \02.13, we then have the following, which is
Theorem \00.10.

\prop{\03.3}  If $Y_1,\dots,Y_q$ weakly intersect properly, then
they have the Autissier property.
\endit

\demo{Proof}  First, note that if
$P\in X\setminus\bigcup_{i\ne j}(Y_i\cap Y_j)$; i.e., if $P\in X$ lies in
at most one of the $Y_i$, then the Autissier property holds trivially at $P$,
because (\03.2.1) is trivial when $r\le 1$.

For all $P\in\bigcup_{i\ne j}(Y_i\cap Y_j)$, we then have that $Y_1,\dots,Y_q$
intersect properly at $P$; therefore they have the Autissier property at $P$
by Proposition \02.13.\qed
\enddemo

(Of course, if $Y_1,\dots,Y_q$ intersect properly, then the first paragraph
of the above proof is unnecessary.)

\beginsection{\04}{Filtrations and Convexity}

This section summarizes the core of Autissier's argument in \citep{a},
as adapted for working with subschemes.

Throughout this section, we fix a complete variety $X$ over $k$ and
proper closed subschemes $Y_1,\dots,Y_q$ of $X$.
Let $\Cal I_1,\dots,\Cal I_q$ be the sheaves of ideals in $\Cal O_X$,
corresponding to $Y_1,\dots,Y_q$, respectively.

We start with some definitions.

\defn{\04.1}  Let $\square$ and $N(\bold t,x)$ be as in Proposition \02.14.
\roster
\myitem a.  Let $N$ be a saturated subset of $\Bbb N^q$.  Then
$$\Cal J_X(N) = \sum_{\bold b\in N} \Cal I_1^{b_1}\dotsm\Cal I_q^{b_q}\;.
  \tag\04.1.1$$
This is a coherent ideal sheaf in $\Cal O_X$.
\myitem b.  For each $\bold t\in\square$ and all $x\in\Bbb R_{\ge0}$, let
$$\Cal J_X(\bold t,x) = \Cal J_X(N(\bold t,x))
  = \sum_{\bold b\in N(\bold t,x)} \Cal I_1^{b_1}\dotsm\Cal I_q^{b_q}\;.
  \tag\04.1.2$$
\myitem c.  Fix a line sheaf $\Cal L$ on $X$, and let $\bold t$ and $x$
be as above.  Then we let
$$\Cal F(\bold t)_x = \Cal F_{\Cal L}(\bold t)_x
  = H^0(X,\Cal L\otimes\Cal J_X(\bold t,x))\;.\tag\04.1.3$$
Then $(\Cal F(\bold t)_x)_{x\in\Bbb R_{\ge0}}$ is a descending filtration
of $H^0(X,\Cal L)$ that satisfies $\Cal F(\bold t)_x=0$ for all $x\gg0$.
\myitem d.  Finally, for all $\bold t\in\square$ we let
$$F(\bold t) = F_{\Cal L}(\bold t)
  = \frac1{h^0(X,\Cal L)} \int_0^\infty \bigl(\dim\Cal F(\bold t)_x\bigr)\,dx\;.
  \tag\04.1.4$$
\endroster
\endit

In terms of this definition, Proposition \02.14 gives the following.

\lemma{\04.3}  Assume that $Y_1,\dots,Y_q$ have the Autissier property,
and let $\Cal L$ be a line sheaf on $X$.  Let $\square$ and $N(\bold t,x)$
be as in Proposition \02.14.  Then
$$\Cal F(\bold t)_x\cap\Cal F(\bold u)_y
  \subseteq \Cal F(\lambda\bold t+(1-\lambda)\bold u)_{\lambda x+(1-\lambda)y}
  \tag\04.3.1$$
for all $\bold t,\bold u\in\square$, all $x,y\in\Bbb R_{\ge0}$, and
all $\lambda\in[0,1]$.
\endit

\demo{Proof}  Let $\bold t,\bold u,x,y,\lambda$ be as above.
By Proposition \02.14 (applied at all $P\in X$),
$$\Cal J_X(\bold t,x) \cap \Cal J_X(\bold u,y)
  \subseteq \Cal J_X(\lambda\bold t+(1-\lambda)\bold u,
    \lambda x+(1-\lambda)y)\;.\tag\04.3.2$$
This remains true after tensoring with $\Cal L$, and (\04.3.1) then follows
because the global section functor is left exact.\qed
\enddemo

We then have the following concavity theorem of
\citet{a, Th\'eor\`eme 3.6} (see also \citep{rv1, Prop.~6.7}).

\thm{\04.4}  Let $\Cal F(\bold t)_x$ ($\bold t\in\square$, $x\in\Bbb R_{\ge0}$)
and $F\:\square\to\Bbb R$ be as in Definition \04.1.
Let $\beta_1,\dots,\beta_q\in\Bbb R_{>0}$.  If (\04.3.1) holds, then
the inequality
$$F(\bold t)
  \ge \min_{1\le i\le q} \left(\frac1{\beta_i}\sum_{m=1}^\infty
    \frac{h^0(X,\Cal L\otimes\Cal I_i^m)}{h^0(X,\Cal L)}\right)\tag\04.4.1$$
holds for all $\bold t\in\square$ for which $\sum \beta_it_i=1$.
\endit

\demo{Proof}  See \citep{rv1, Prop.~6.7}.\qed
\enddemo

The results of this section can then be summarized as follows.

\thm{\04.5}  If $Y_1,\dots,Y_q$ have the Autissier property
and if $\beta_1,\dots,\beta_q\in\Bbb R_{>0}$, then (\04.4.1) holds.
\endit

This provides a slight strengthening of the ``General Theorem'' of \citep{rv1}:
in that theorem, the divisors $D_i$ were assumed to be Cartier, but this
condition has been relaxed so that they only need to be Cartier
at points where they meet other divisors in the collection.

\beginsection{\05}{Ideal Sheaves and B-divisors}

The remainder of the proof in \citep{rv1} involves Prop.~4.18 of that paper,
so it is necessary to interpret things such as
$H^0(X,\Cal L\otimes\Cal I_1^{b_1}\dots\Cal I_q^{b_q})$
in terms of Cartier b-divisors.  This is quite easy, because ideal sheaves
are special cases of b-divisors.  That is the topic of this section.

Throughout this section, $X$ is a complete variety over a field $k$,
unless otherwise specified.

We briefly recall that a {\bc model} of $X$
is a proper birational morphism $\pi\:W\to X$ of varieties over $k$,
and a {\bc Cartier b-divisor} $\bold D$ on $X$ is an equivalence class
of pairs $(W,D)=(\pi\:W\to X,D)$, where $\pi\:W\to X$ is a model of $X$
and $D$ is a Cartier divisor on $W$; here pairs $(W_1,D_1)$ and $(W_2,D_2)$
are said to be equivalent if there exist a model $W_3$ of $X$ and morphisms
$f_i\:W_3\to W_i$ over $X$ for $i=1,2$ such that $f_1^{*}D_1=f_2^{*}D_2$.
For more details and basic properties, see \citep{rv1, \S\,4}.

We start with some basic results about spaces of global sections of
line sheaves on projective varieties.  The first result is a general result
on growth of cohomology groups, and is essentially due to
the \citet{stks, Lemma 0BEM}.
The latter lemma says that the Euler characteristic
of the sheaves $\Cal F\otimes\Cal L_1^{n_1}\otimes\dots\otimes\Cal L_r^{n_r}$
is a numerical polynomial in $n_1,\dots,n_r$ of a certain degree.
Although the lemma below gives instead a bound on the dimensions of the
cohomology groups of these sheaves, the method of proof is the same.
These upper bounds will only be needed for $h^0$, but we will prove
the general case as it is no more difficult.

\lemma{\05.1}  Let $X$ be a proper scheme over a field $k$,
let $\Cal F$ be a coherent sheaf on $X$, let $d=\dim\Supp\Cal F$, and
let $\Cal L_1,\dots,\Cal L_r$ be line sheaves on $X$.  Then
$$h^i(X,\Cal F\otimes\Cal L_1^{n_1}\otimes\dots\otimes\Cal L_r^{n_r})
  \le O(|\bold n|^d+1)$$
for all $\bold n=(n_1+\dots+n_r)\in\Bbb N^r$ and all $i$, where the implicit
constant depends only on $X$, $k$, $\Cal F$, and $\Cal L_1,\dots,\Cal L_r$.
\endit

\demo{Proof}  We give a sketch of this proof, following \citep{stks},
including all places where the proofs differ.

For typographical simplicity, we let $\Cal L^{\bold n}$ denote
$\Cal L_1^{n_1}\otimes\dots\otimes\Cal L_r^{n_r}$ (multiindex notation)
for all $\bold n\in\Bbb N^r$.

The proof is by induction on $d$.  The base case $d=0$
(including also $\Supp\Cal F=\emptyset$) is trivial.

First, if $\Cal F$ contains embedded points, then by \citep{stks, Lemma 02OL}
there is a short exact sequence
$$0 @>>> \Cal K @>>> \Cal F @>>> \Cal F' @>>> 0\tag\05.1.1$$
of coherent sheaves
such that $\dim\Supp\Cal K<d$ and $\Cal F'$ has no embedded points.
It remains exact after tensoring with $\Cal L^{\bold n}$, so by the long exact
sequence in cohomology and the inductive hypothesis we have
$$\left|h^i(X,\Cal F\otimes\Cal L^{\bold n})
    - h^i(X,\Cal F'\otimes\Cal L^{\bold n})\right|
  \le O(|\bold n|^{d-1}+1)\;.$$
Therefore it suffices to prove the lemma when $\Cal F$ has no embedded points.

We may replace $X$ with $\Supp\Cal F$ (this does not change the cohomologies),
so we may assume that $\dim X=d$ and that $X$ has no embedded points.
In this situation, by \citep{stks, Lemmas~02OZ and~02P2}, there exist
a coherent ideal sheaf $\Cal I$ on $X$ and short exact sequences
$$0 @>>> \Cal I\Cal F @>>> \Cal F @>>> \Cal Q @>>> 0
  \qquad\text{and}\qquad
  0 @>>> \Cal I\Cal F @>>> \Cal F\otimes\Cal L_1 @>>> \Cal Q' @>>> 0$$
such that $\dim\Supp\Cal Q<d$ and $\dim\Supp\Cal Q'<d$.  Again tensoring
with $\Cal L^{\bold n}$ and applying the long exact sequence and
the inductive hypothesis, we have
$$\left|h^i(X,\Cal F\otimes\Cal L^{\bold n})
    - h^i(X,\Cal I\Cal F\otimes\Cal L^{\bold n})\right|
  \le O(|\bold n|^{d-1}+1)$$
and
$$\left|h^i(X,\Cal F\otimes\Cal L^{\bold n}\otimes\Cal L_1)
    - h^i(X,\Cal I\Cal F\otimes\Cal L^{\bold n})\right|
  \le O(|\bold n|^{d-1}+1)$$
for all $\bold n$ and all $i$.  Combining these inequalities, and using
a symmetrical argument, we obtain
$$\left|h^i(X,\Cal F\otimes\Cal L^{\bold n}\otimes\Cal L_j)
    - h^i(X,\Cal F\otimes\Cal L^{\bold n})\right|
  \le O(|\bold n|^{d-1}+1)$$
for all $\bold n$, all $i$, and all $j=1,\dots,r$.

Applying this inequality $|\bold n|$ times then gives
$$\left|h^i(X,\Cal F\otimes\Cal L^{\bold n}) - h^i(X,\Cal F)\right|
  \le O(|\bold n|^d+1)\;,$$
and the result follows.\qed
\enddemo

The following lemma applies this to give bounds more directly applicable
to the current situation.

\lemma{\05.2}  Let $\pi\:W'\to W$ be a proper birational morphism of complete
varieties over a field $k$.
\roster
\myitem a.  Assume that $W$ is normal.  Then $\pi_{*}\pi^{*}\Cal L\cong\Cal L$
for all line sheaves $\Cal L$ on $W$, and the natural map
$H^0(W,\Cal L)\to H^0(W',\pi^{*}\Cal L)$ is an isomorphism.
\myitem b.  For general $W$, the coherent sheaf
$\Cal F=\pi_{*}\Cal O_{W'}/\Cal O_W$ on $W$ is supported on a proper subset
of $W$, and
$$0 \le h^0(W',\pi^{*}\Cal L) - h^0(W,\Cal L)
  \le h^0(W,\Cal F\otimes\Cal L)\tag\05.2.1$$
for all line sheaves $\Cal L$ on $W$.
\myitem c.  Let $\Cal L$ be a line sheaf on $W$, let $D$ be a Cartier divisor
on $W$, and let $d=\dim W$.  Then
$$0 \le h^0(W',\pi^{*}\Cal L^N(-m\pi^{*}D)) - h^0(W,\Cal L^N(-mD))
  \le O((N+m)^{d-1})\tag\05.2.2$$
for all $N\in\Bbb Z_{>0}$ and all $m\in\Bbb N$, where the implicit constant
depends on $\pi$, $k$, $\Cal L$, and $D$, but not on $N$ or $m$.
\myitem d.  Under the same conditions as (c),
$$0 \le \sum_{m=1}^\infty h^0(W',\pi^{*}\Cal L^N(-m\pi^{*}D))
    - \sum_{m=1}^\infty h^0(W,\Cal L^N(-mD))
  \le O(N^d)\;.\tag\05.2.3$$
\endroster
\endit

\demo{Proof}  For part (a), we first note that $\pi_{*}\Cal O_{W'}=\Cal O_W$
(as subsheaves on the constant sheaves of the function field $K(W')\cong K(W)$)
by \citep{ha77, II Prop.~6.3A} and the fact that $W$ is normal.
Therefore the projection formula gives $\pi_{*}\pi^{*}\Cal L\cong\Cal L$,
and taking global sections gives $H^0(W',\pi^{*}\Cal L)\cong H^0(W,\Cal L)$.

For (b), we have an exact sequence
$$ 0 @>>> \Cal O_W @>>> \pi_{*}\Cal O_{W'} @>>> \Cal F @>>> 0$$
of sheaves on $W$, where $\Cal F$ is supported on a proper subset of $W$.
Tensoring each term with $\Cal L$ and taking global sections then gives
an exact sequence
$$0 @>>> H^0(W,\Cal L) @>>> H^0(W',\pi^{*}\Cal L)
  @>>> H^0(W,\Cal F\otimes\Cal L)\;,\tag\05.2.4$$
which gives (\05.2.1).

By (b), part (c) is a matter of showing that
$$h^0(W,\Cal F\otimes\Cal L^N(-mD)) \le O((N+m)^{d-1})$$
for all $N$ and $m$.  This is immediate from Lemma \05.1
with $\Cal L_1=\Cal L$ and $\Cal L_2=\Cal O(-D)$, since
$\dim\Supp\Cal F\le d-1$.

For (d), the lower bound holds (termwise) by the first part of (\05.2.2).

For the upper bound, we first note that there is a constant $c$
(independent of $N$ and $m$) such that the summands in (\05.2.3)
vanish for all $m>cN$.  Indeed, let $\pi''\:W''\to W'$ be a projective model
of $W$ that dominates $W'$ and let $A$ be an ample divisor on $W''$;
then it suffices to take
$c \ge
 (\pi^{\prime\prime*}\Cal L\idot A^{d-1})/(\pi^{\prime\prime*}D\idot A^{d-1})$,
where in this case $A^{d-1}$ is meant in the sense of intersection theory.

The sums then have $O(N)$ nonzero terms with $m\le O(N)$, so the upper bound
follows from (\05.2.2).\qed
\enddemo

\defn{\05.3}  Let $\Cal L$ be a line sheaf on $X$ and let $\bold D$ be
an effective Cartier b-divisor on $X$.  Then
$$H^0_\bir(X,\Cal L(-\bold D)) = H^0(W,\pi^{*}\Cal L(-D))\;,$$
where $\pi\:W\to X$ is any normal model of $X$ on which $\bold D$
is represented by a Cartier divisor $D$.  This is independent
of the choice of $W$ by Lemma \05.2a.

Also (as usual)
$$h^0_\bir(X,\Cal L(-\bold D))=\dim_k H^0_\bir(X,\Cal L(-\bold D))\;.$$

When $\bold D=0$, these are also denoted $H^0_\bir(X,\Cal L)$
and $h^0_\bir(X,\Cal L)$, respectively.
\endit

The subscript ``bir'' is needed because $H^0_\bir(X,\Cal L)$ may differ
from $H^0(X,\Cal L)$ if $X$ is not normal.

\lemma{\05.4}  Let $\Cal L$ be a line sheaf on $X$, let $D$ be a
nonzero effective Cartier divisor on $X$, and let $d=\dim X$.  Then
$$h^0_\bir(X,\Cal L^N) = h^0(X,\Cal L^N) + O(N^{d-1})\tag\05.4.1$$
and
$$\sum_{m=1}^\infty h^0_\bir(X,\Cal L^N(-mD))
  = \sum_{m=1}^\infty h^0(X,\Cal L^N(-mD)) + O(N^d)\tag\05.4.2$$
as $N\to\infty$, where the implicit constants depend only on $\Cal L$ and $D$.
In particular, if $\Cal L$ is big, then
$$\beta(\Cal L, D)
  = \liminf_{N\to\infty} \frac{\sum_{m=1}^\infty h^0_\bir(X,\Cal L^N(-mD))}
    {Nh^0_\bir(X,\Cal L^N)}\;.\tag\05.4.3$$
\endit

\demo{Proof}  First of all, by Lemma \05.2a, $h^0_\bir(X,\Cal L^N(-mD))$
for all $N,m\in\Bbb N$ can be computed on a fixed normal model $W$ of $X$,
independent of $N$ and $m$.

Then (\05.4.1) is immediate from Lemma \05.2c.

For (\05.4.2), note that $h^0_\bir(X,\Cal L^N(-mD))=h^0(W,\pi^{*}\Cal L^N(-mD))$
for any normal model $\pi\:W\to X$.  Then (\05.4.2) is immediate from
Lemma \05.2d applied to any such model $\pi$.

Finally, since $\Cal L$ is big, (\05.4.3) follows easily from
(\05.4.1) and (\05.4.2).\qed
\enddemo

Therefore we may extend Definition \00.1 as follows.

\defn{\05.5}  Let $\Cal L$ be a big line sheaf on $X$ and let $\bold D$
be a nonzero effective Cartier b-divisor on $X$.  Then
$$\beta(\Cal L, \bold D)
  = \liminf_{N\to\infty}
    \frac{\sum_{m=1}^\infty h^0_\bir(X,\Cal L^N(-m\bold D))}
      {Nh^0_\bir(X,\Cal L^N)}\;.\tag\05.5.1$$
\endit

\remk{\05.6}  As noted in \citep{rv1} (following Def.~1.9),
the above $\liminf$ is actually a limit when $\Cal L$ is big and $D$ is
a Cartier divisor.  A detailed proof is given in Section \010, including
the case when $\bold D$ is a b-divisor.
\endit

Now we consider b-divisors associated to proper closed subschemes.

\defn{\05.7}  Let $Y$ be a proper closed subscheme of $X$, and let $\Cal I$
be the corresponding ideal sheaf.  Let $\pi\:W\to X$ be the blow-up of $X$
along $\Cal I$, and let $E$ be the exceptional divisor of $\pi$ (so that
$\Cal O(E)=\Cal O(-1)$ for the blowing-up).  Then the Cartier b-divisor
$\bold Y$ associated to $Y$ is the b-divisor represented by $E$ on $W$.
\endit

Next we compare the relevant spaces of global sections.

\lemma{\05.8}  Let $\Cal L$ be a line sheaf on $X$.
\roster
\myitem a.  Let $Y$, $\Cal I$, $\pi\:W\to X$, and $E$ be as in
Definition \05.7.  Let $m\in\Bbb N$.  Then the restriction to $\Cal I^m$
of the natural map $\Cal O_X\hookrightarrow\pi_{*}\Cal O_W$ gives a map
$$\Cal I^m\hookrightarrow \pi_{*}(\Cal O_W(-mE))\;.\tag\05.8.1$$
This gives an injection
$$H^0(X,\Cal L\otimes\Cal I^m) \hookrightarrow H^0(W,\pi^{*}\Cal L(-mE))\;.
  \tag\05.8.2$$
\myitem b.  The map (\05.8.2) is an isomorphism for all sufficiently large $m$,
independent of $\Cal L$.
\myitem c.  For each $i=1,\dots,q$ let $Y_i$ be a proper closed subscheme
of $X$, and let $\bold Y_i$ and $\Cal I_i$ be the corresponding
Cartier b-divisor and ideal sheaf on $X$, respectively.
Then, for all $n_1,\dots,n_q\in\Bbb N$, there is a canonical injection
$$H^0(X,\Cal L\otimes\Cal I_1^{n_1}\dotsm\Cal I_q^{n_q})
  \hookrightarrow
  H^0_\bir(X,\Cal L(-n_1\bold Y_1-\dots-n_q\bold Y_q))\;,\tag\05.8.3$$
induced by the maps of part (a) for all $i$.
\endroster
\endit

\demo{Proof}  With notation as in part (a), let $U$ be an open subset
of $X$.  Then any local section $s\in\Gamma(U,\Cal I)$ pulls back to a
section of $\pi^{-1}\Cal I\cdot\Cal O_W=\Cal O(1)=\Cal O(-E)$
over $\pi^{-1}(U)$; see \citep{ha77, II Prop.~7.13}.  This gives (\05.8.1).

Tensoring both sides of (\05.8.1) with $\Cal L$ and applying the
projection formula gives an injection
$\Cal L\otimes\Cal I^m\hookrightarrow\pi_{*}(\pi^{*}\Cal L(-mE))$,
which then gives (\05.8.2).

For part (b), it suffices to show that the map (\05.8.1) is surjective
(hence an isomorphism) for all $m\gg0$.  This map can be written
$\Cal I^m\to\pi_{*}\Cal O_W(m)$.  The fact that it is surjective for all
$m\gg0$ is noted at the very end of the proof of \citep{ha77, II Thm.~5.19}.
(This is shown locally over open affines of $X$, but extends to all of $X$
by a compactness argument.)

For part (c), let $\pi\:W\to X$ be any normal model of $X$ that dominates
the blowings-up of $X$ along $Y_i$ for all $i$.
Since $\Cal L\otimes\Cal I_1^{n_1}\dotsm\Cal I_q^{n_q}$ is locally generated
by products of local sections of $\Cal L$ and
of $\Cal I_1^{n_1},\dots,\Cal I_q^{n_q}$, we obtain from (\05.8.1) an injection
$\Cal L\otimes\Cal I_1^{n_1}\dotsm\Cal I_q^{n_q}
  \hookrightarrow \pi_{*}(\pi^{*}\Cal L(-n_1E_1-\dots-n_qE_q))$,
which gives (\05.8.3).\qed
\enddemo

We conclude this section by proving the assertions of Remarks \00.5 and \00.6.

\cor{\05.9}  Let $Y$ be a proper closed subscheme of $X$ and let $\bold Y$
be the corresponding b-divisor.  Let $\Cal L$ be a big line sheaf on $X$.
Then:
\roster
\myitem a.  Recalling Definitions \00.4, \00.5, and \05.5,
$$\beta(\Cal L,Y) = \beta_{\Cal L,Y} = \beta(\Cal L,\bold Y)\;.\tag\05.9.1$$
\myitem b.  If any of these three quantities can be computed by evaluating
the corresponding limits, then all of them can.
\endroster
\endit

\demo{Proof}  Let $\Cal I$, $\pi\:W\to X$, and $E$ be as in Definition \05.7,
and let $d=\dim X$.  For all $m\in\Bbb Z_{>0}$ let $\Cal F_m$ be the cokernel
of the map (\05.8.1); by Lemma \05.8b there is an $m_0$ such that
$\Cal F_m=0$ for all $m>m_0$.  Tensoring the short exact sequence
$0 @>>> \Cal I^m @>>> \pi_{*}\Cal O_W(-mE) @>>> \Cal F_m @>>> 0$
with $\Cal L^N$ and taking global sections then gives
$$\split 0 &\le \sum_{m=1}^\infty h^0(W,\pi^{*}\Cal L^N(-mE))
    - \sum_{m=1}^\infty h^0(X,\Cal L^N\otimes\Cal I^m) \\
  &\le \sum_{m=1}^{m_0} h^0(X,\Cal F_m\otimes\Cal L^N) \\
  &\le O(N^{d-1}+1)\endsplit$$
for all $N>0$, by Lemma \05.2c.  This gives
$$\liminf_{N\to\infty}
    \frac{\sum_{m=1}^\infty h^0(X,\Cal L^N\otimes\Cal I^m)}{Nh^0(X,\Cal L^N)}
  = \liminf_{N\to\infty}
    \frac{\sum_{m=1}^\infty h^0(W,\pi^{*}\Cal L^N(-mE))}{Nh^0(X,\Cal L^N)}\;,
  \tag\05.9.2$$
which is the first equality $\beta(\Cal L,Y) = \beta_{\Cal L,Y}$ of (\05.9.1).

The second equality $\beta_{\Cal L,Y} = \beta(\Cal L,\bold Y)$ is a matter
of showing that
$$\liminf_{N\to\infty}
    \frac{\sum_{m=1}^\infty h^0(W,\pi^{*}\Cal L^N(-mE))}{Nh^0(X,\Cal L^N)}
  = \liminf_{N\to\infty}
    \frac{\sum_{m=1}^\infty h^0_\bir(X,\Cal L^N(-m\bold D))}
      {Nh^0_\bir(X,\Cal L^N)}\tag\05.9.3$$
This is true by (\05.4.2) and (\05.4.1).

Part (b) is immediate from the fact that (\05.9.2) and (\05.9.3) remain true
(for the same reasons) when all instances of $\liminf$ are replaced
by $\limsup$.\qed
\enddemo

\beginsection{\06}{An Inequality of B-divisors}

This section continues with the proof of Theorem \00.11, by applying
the method of \citet{a, \S\,4} as adapted in \citep{rv1},
leading up to an inequality of $\Bbb R$\snug-Cartier b-divisors (Lemma \06.9).
This closely follows the proof in \citep{rv1, \S\,6}, but we simplify it here
by eliminating the sets $\Sigma$ and $\triangle_\sigma$ (see Remark \06.10).

We start with some notation.  Let $X$ and $Y_1,\dots,Y_q$ be as in
the statement of Theorem \00.11, and let $\beta_1,\dots,\beta_q\in\Bbb R_{>0}$.
Let $b$ and $N$ be large positive integers, to be chosen later
(Proposition \07.4).

Let
$$\triangle = \{\bold t\in\Bbb R_{\ge0}^q:t_1+\dots+t_q=1\}\;.$$
Recalling that $b\in\Bbb Z_{>0}$, let
$$\triangle_b
  = \left\{\bold a\in\prod_{i=1}^q \beta_i^{-1}\Bbb N
    : \sum\beta_ia_i=b\right\}\;,$$
so that
$b^{-1}\triangle_b$ is a finite discrete subset of $\triangle$.

Recall from Sections \02 and \04 that
$\square=\Bbb R_{\ge0}^q\setminus\{\bold 0\}$ and that
$$N(\bold t,x)=\{\bold b\in\Bbb N^q:\sum t_ib_i\ge x\}\;,
  \qquad\bold t\in\square\,,\; x\in\Bbb R_{\ge0}\;.$$
Let $F=F_{\Cal L^N}\:\square\to\Bbb R$ be the function of Definition \04.1.
Write $b^{-1}\bold a=\bold a/b$ for all $\bold a\in\triangle_b$,
and recall that
$$F(\bold a/b)
  = \int_0^\infty
    \frac{\dim\Cal F_{\Cal L^N}(\bold a/b)_x}{h^0(X,\Cal L^N)}\,dx\;,$$
where $(\Cal F(\bold a/b))_x = (\Cal F_{\Cal L^N}(\bold a/b))_x$ is
the filtration of $H^0(X,\Cal L^N)$ given by
$$\Cal F(\bold a/b)_x = H^0(X,\Cal L^N\otimes\Cal J_X(\bold a/b,x))$$
and
$$\Cal J_X(\bold a/b,x)
  = \sum_{\bold b\in N(\bold a/b,x)} \Cal I_1^{b_1}\dotsm\Cal I_q^{b_q}\;.$$

For all $\bold a\in\triangle_b$ and $x\in\Bbb R_{\ge0}$
let $K=K(\bold a/b,x)$ be the set of minimal elements in $N(\bold a/b,x)$.
Then
$$\Cal J_X(\bold a/b,x)
  = \sum_{\bold b\in K} \Cal I_1^{b_1}\dotsm\Cal I_q^{b_q}\;,\tag\06.1$$
and this is a finite sum since $K$ is a finite set.

Following \citet{rv1, \S\,6}, for all $\bold a\in\triangle_b$
and all $s\in H^0(X,\Cal L^N)\setminus\{0\}$ we define
$$\mu_{\bold a/b}(s) = \sup\{x:\Cal F(\bold a/b)_x\ni s\}\;.\tag\06.2$$

\lemma{\06.3}  Let $\bold a\in\triangle_b$
and $s\in H^0(X,\Cal L^N)\setminus\{0\}$.  Let $\mu=\mu_{\bold a/b}(s)$.  Then
$$s \in
  H^0\left(X,\sum_{\bold b\in K(\bold a/b,\mu)}
     \Cal L^N\otimes\Cal I_1^{b_1}\dotsm\Cal I_q^{b_q}\right)\;.\tag\06.3.1$$
\endit

\demo{Proof}  The union $\bigcup_{x\in[0,\mu]} K(\bold a/b,x)$ is finite,
and each $\bold b$ in this union occurs in the sum (\06.1)
for a closed set of $x$.  Therefore the supremum in (\06.2) is actually
a maximum.  In particular, $s\in\Cal F(\bold a/b)_\mu$, and this gives
(\06.3.1).\qed
\enddemo

\remk{\06.4}  Since the injection (\05.8.3) is not necessarily bijective,
it is important in this section to carefully distinguish between objects
defined on $X$ (non-birational objects) and the birational objects
defined in Section \05.  So far in this Section \06, everything
has been non-birational.  This will now change.
\endit

\cor{\06.5}  Let $\bold a$, $s$, and $\mu$ be as in Lemma \06.3, and let
$K=K(\bold a/b,s)$.  Let $\bold Y_1,\dots,\bold Y_q$ be the b-divisors on $X$
corresponding to $Y_1,\dots,Y_q$, respectively.  Then
$$(s) \ge \bigwedge_{\bold b\in K} \sum_{i=1}^q b_i\bold Y_i\;.\tag\06.5.1$$
\endit

\demo{Proof}  Let $\pi\:W\to X$ be a model of $X$ on which all $\bold Y_i$
are represented by Cartier divisors $D_i$.  Then, by (\06.3.1),
$s$ is a global section of the subsheaf of $\pi^{*}\Cal L^N$ generated by
the set $\{\pi^{*}\Cal L^N(-b_1D_1-\dots-b_qD_q):\bold b\in K\}$.
By \citep{rv1, Prop.~4.18}, since this set is finite, we have
$$(\pi^{*}s) \ge \bigwedge_{\bold b\in K}(b_1D_1+\dots+b_qD_q)\;.$$
This gives (\06.5.1).\qed
\enddemo

\defn{\06.6}  Let $\Cal F=(\Cal F_x)_{x\in\Bbb R_{\ge0}}$ be a filtration
of a finite dimensional vector space $V$, and let $\Cal B$ be a basis of $V$.
Then $\Cal B$ is {\bc adapted} to $\Cal F$ if $\Cal B\cap\Cal F_x$ is a basis
of $\Cal F_x$ for all $x$.
\endit

\defn{\06.7}  Let $\Cal B$ be a basis of $H^0(X,\Cal L^N)$.  Then
$$(\Cal B) = \sum_{s\in\Cal B} (s)\;.$$
\endit

\remk{\06.8}  At this point we start using $\Bbb R$\snug-Cartier b-divisors.
These are basically finite formal linear combinations of Cartier b-divisors
with real coefficients.
An $\Bbb R$\snug-Cartier b-divisor is said to be {\bc effective} if it is
a finite linear combination of effective Cartier b-divisors with positive
(real) coefficients.  For more details on $\Bbb R$\snug-Cartier divisors
and $\Bbb R$\snug-Cartier b-divisors, see \citep{rv2, \S\,2}.
\endit

\lemma{\06.9}  Assume that $\operatorname{char} k=0$.
For each $\bold a\in\triangle_b$ let $\Cal B_{\bold a}$ be a
basis of $H^0(X,\Cal L^N)$ adapted to the filtration $\Cal F(\bold a/b)$.
Then
$$\bigvee_{\bold a\in\triangle_b} (\Cal B_{\bold a})
  \ge \frac b{b+q}
    \left(\min_{1\le i\le q} \sum_{m=1}^\infty
      \frac{h^0(X,\Cal L^N\otimes\Cal I_i^m)}{\beta_i}\right)
    \sum_{i=1}^q \beta_i\bold Y_i\;.\tag\06.9.1$$
\endit

\demo{Proof}  Let $\bold D'$ be the left-hand side of (\06.9.1),
and let $\pi\:W\to X$ be a model of $X$ on which $\bold D'$ and
$\bold Y_1,\dots,\bold Y_q$ are represented by Cartier divisors $D'$
and $D_1,\dots,D_q$, respectively.  We also assume that $W$ is nonsingular.

Let $E$ be a prime divisor on $W$.  Let $\nu'$, $\nu_{\bold a}$
for all $\bold a\in\triangle_b$, and $\nu_1,\dots,\nu_q$ be the
multiplicities of $E$ in $D'$, $(\Cal B_{\bold a})$ for all $\bold a$,
and $D_1,\dots,D_q$, respectively.  Let $\nu=\sum\beta_i\nu_i$.

We claim that there is an $\bold a\in\triangle_b$ (depending on $E$)
such that
$$\nu_{\bold a} \ge \frac b{b+q} h^0(X,\Cal L^N)F(\bold a/b)\nu\;.\tag\06.9.2$$

Since the divisor $(\Cal B_{\bold a})$ is effective for all $\bold a$
(and $\triangle_b$ is nonempty), the claim is trivial if $\nu=0$,
so we assume that $\nu>0$.

Let
$$t_i = \frac{\nu_i}{\nu}\;,\qquad i=1,\dots,q\;.\tag\06.9.3$$
Since $\sum\beta_i\nu_i=\nu$, we have $\sum \beta_it_i=1$ and therefore
$b \le \sum\lfloor(b+q)\beta_it_i\rfloor \le b+q$.
Therefore we may choose $\bold a\in\triangle_b$ such that
$$a_i \le (b+q)t_i\;,\qquad i=1,\dots,q\;.\tag\06.9.4$$

Let $s\in\Cal B_{\bold a}$, and let $\nu_s$ be the multiplicity of $E$
in the divisor $(\pi^{*}s)$.  Let $K=K(\bold a/b,\mu_{\bold a/b}(s))$.
By (\06.5.1), (\06.9.3), (\06.9.4),
and the fact that $\sum a_ib_i\ge b\mu_{\bold a/b}(s)$
for all $\bold b\in K\subseteq N(\bold a/b,\mu_{\bold a/b}(s))$,
$$\frac{\nu_s}{\nu} \ge \frac1\nu \min_{\bold b\in K} \sum_{i=1}^q b_i\nu_i
  = \min_{\bold b\in K} \sum_{i=1}^q b_it_i
  \ge \min_{\bold b\in K} \sum_{i=1}^q \frac{a_ib_i}{b+q}
  \ge \frac{b}{b+q}\mu_{\bold a/b}(s)\;.\tag\06.9.5$$

Since $\Cal B_{\bold a}$ is adapted to the filtration $\Cal F(\bold a/b)$,
we have
$$h^0(X,\Cal L^N)F(\bold a/b) = \int_0^\infty \dim\Cal F(\bold a/b)_x\,dx
  = \sum_{s\in\Cal B_{\bold a}} \mu_{\bold a/b}(s)$$
(see \citep{rv1, Remark~6.6}).  Combining this with (\06.9.5) and
the fact that $\nu_{\bold a}=\sum_{s\in\Cal B_{\bold a}}\nu_s$ then gives
(\06.9.2).

Since $Y_1,\dots,Y_n$ have the Autissier property, Theorem \04.5 gives
$$h^0(X,\Cal L^N)F(\bold a/b)
  \ge \min_{1\le i\le q} \sum_{m=1}^\infty
    \frac{h^0(X,\Cal L\otimes\Cal I_i^m)}{\beta_i}\;.$$
Therefore, by (\06.9.2) and the definition of $\nu$, we have
$$\nu'
  \ge \frac b{b+q}
    \left(\min_{1\le i\le q} \sum_{m=1}^\infty
      \frac{h^0(X,\Cal L^N\otimes\Cal I_i^m)}{\beta_i}\right)
    \sum_{i=1}^q \beta_i\nu_i\;.$$
We conclude that the difference of the two sides of (\06.9.1) is represented
on $W$ by a finite sum
of effective Cartier divisors with nonnegative real coefficients
(these divisors are the finitely many prime divisors $E$ occurring in
$\Supp D'$, and they are Cartier because $W$ is nonsingular).
This proves (\06.9.1).\qed
\enddemo

\remk{\06.10}  As noted in the introductory paragraph of this section,
we have simplified the argument somewhat by eliminating the dependence
on subsets $\sigma\subseteq\{1,\dots,q\}$.  It would be easy to put this
dependence back (by Remark \02.10 it is still true that at most $\dim X$
of the $Y_i$ can pass through any point of $X$).  With this change,
the fraction $b/(b+q)$ in Lemma \06.9 can be replaced by $b/(b+\dim X)$,
as in \citep{rv1}.
\endit

\beginsection{\07}{A Birational Nevanlinna Constant for B-divisors}

In this section we introduce the birational Nevanlinna constant of \citet{rv1},
as modified to use $\Bbb R$\snug-Cartier b-divisors, and prove the bound
(\07.4.2), which corresponds to the penultimate step in the proof of
the Main Theorem of \citep{rv1}.

We start with the following definition, which is \citep{rv2, Def.~1.1},
except that $\bold D$ is allowed to be an $\Bbb R$\snug-Cartier b-divisor
instead of an $\Bbb R$\snug-Cartier divisor.

\defn{\07.1}  Let $X$ be a complete variety, let $\Cal L$ be a line sheaf
on $X$, and let $\bold D$ be an effective $\Bbb R$\snug-Cartier b-divisor
on $X$.  Then
$$\Nevbir(\Cal L,\bold D) = \inf_{N,V,\mu} \frac{\dim V}{\mu}\;,\tag\07.1.1$$
where the infimum passes over all triples $(N,V,\mu)$ such that
$N\in\Bbb Z_{>0}$, $V$ is a linear subspace of $H^0(X,\Cal L^N)$
with $\dim V>1$, and $\mu\in\Bbb R_{>0}$, with the following property.
There exists a model $\pi\:W\to X$ such that the following condition holds.
For all $Q\in W$ there is a basis $\Cal B$ of $V$ such that
$$(\Cal B) \ge \mu N\bold D\tag\07.1.2$$
in a Zariski-open neighborhood $U$ of $Q$,
relative to the cone of effective $\Bbb R$\snug-Cartier b-divisors on $U$.
If there are no such triples $(N,V,\mu)$, then $\Nevbir(\Cal L,D)$
is defined to be $+\infty$.
\endit

As in \citep{rv1, Cor.~4.17}, we have the following alternative
characterization.

\prop{\07.2}  Let $X$, $\Cal L$, and $\bold D$ be as in Definition \07.1.  Then
$$\Nevbir(\Cal L,\bold D) = \inf_{N,V,\mu} \frac{\dim V}{\mu}\;,$$
where the infimum passes over all triples $(N,V,\mu)$ such that
$N\in\Bbb Z_{>0}$, $V$ is a linear subspace of $H^0(X,\Cal L^N)$
with $\dim V>1$, and $\mu\in\Bbb R_{>0}$, with the following property.
There is a finite list $\Cal B_1,\dots\Cal B_\ell$ of bases
of $H^0(X,\Cal L^N)$ such that
$$\bigvee_{i=1}^\ell (\Cal B_i) \ge \mu N\bold D$$
(with the same convention if there are no such triples).
\endit

\remk{\07.3}  One could make Definition ``fully birational'' by allowing
$\Cal L$ to be a b-line-sheaf.  Here a {\bc b-line-sheaf} is an element
of $\varinjlim\Pic W$, where the direct limit is over all models $W$ of $X$.
However, this would basically amount to replacing $X$ with some model $W$ on
which the hypothetical b-line-sheaf lies in $\Pic W$, so nothing new would
be introduced.
\endit

With Definition \07.1, we have:

\prop{\07.4}  Let $k$, $X$, $\Cal L$, and $Y_1,\dots,Y_q$ be as in
the statement of Theorem \00.11.  For each $i=1,\dots,q$ let $\bold Y_i$
be the Cartier b-divisor on $X$ corresponding to $Y_i$.  Let
$\beta_1,\dots,\beta_q\in\Bbb R_{>0}$, and let $\bold D$ be the
effective $\Bbb R$\snug-Cartier b-divisor
$\beta_1\bold Y_1+\dots+\beta_q\bold Y_q$.  Then
$$\split \Nevbir(\Cal L,\bold D)
  &\le Nh^0(X,\Cal L^N)\fracwithdelims(){b+q}{b}
    \left(\min_{1\le i\le q} \sum_{m=1}^\infty
      \frac{h^0(X,\Cal L^N\otimes\Cal I_i^m)}{\beta_i}\right)^{-1} \\
  &= \fracwithdelims(){b+q}{b}
    \left(\min_{1\le i\le q} \frac1{\beta_i} \sum_{m=1}^\infty
      \frac{h^0(X,\Cal L^N\otimes\Cal I_i^m)}{Nh^0(X,\Cal L^N)}\right)^{-1}
  \endsplit\tag\07.4.1$$
for all $b\in\Bbb Z_{>0}$ and all $N\in\Bbb Z_{>0}$ such that
$H^0(X,\Cal L^N\otimes\Cal I_i)\ne0$ for all $i$.

In particular, if $\beta_i=\beta(\Cal L,Y_i)$ for all $i$, then
$$\Nevbir(\Cal L,\bold D) \le 1\;.\tag\07.4.2$$
\endit

\demo{Proof}  The inequality (\07.4.1) follows from Lemma \06.9 and
Proposition \07.2, with $V=\allowmathbreak H^0(X,\Cal L^N)$.
For the second inequality, we can omit those $Y_i$ for which
$\beta(\Cal L,Y_i)=0$; then (\07.4.2) follows from (\07.4.1) and (\00.4.1)
by taking $b$ and $N$ sufficiently large.\qed
\enddemo

\beginsection{\08}{Conclusion of the Main Proof}

This section gives the last step of the proof of Theorems \00.11 and \00.9.
This relies on Theorem \08.2, which generalizes \citep{rv1, Thms.~1.4 and~1.5}.
(Theorem \08.2 may be regarded as the Second Main Theorem for
birational Nevanlinna constants.)

We start by introducing Weil functions, and the resulting proximity and
counting functions, for $\Bbb R$\snug-Cartier b-divisors $\bold D$.
Basically, this involves lifting to a model on which the b-divisor in question
is represented by an $\Bbb R$\snug-Cartier divisor, and
using $\Bbb R$\snug-linearity.
This model may not be isomorphic to $X$ over all of $X\setminus\Supp\bold D$
(here $\Supp\bold D$ is defined as $\pi(\Supp D)$, where $\pi\:W\to X$ is a
model of $X$ on which $\bold D$ is represented by a Cartier divisor $D$).
Therefore the domain of the Weil function may be smaller than
one would otherwise expect.

For more details on Weil functions, see \citep{rv1, 2.3},
\citep{la_fdg, Ch.~10}, or \citep{voj11, \S\,8}.

\defn{\08.1}  Let $\bold D$ be an $\Bbb R$\snug-Cartier b-divisor on $X$.
Let $\pi\:W\to X$ be a model of $X$ on which $\bold D$ is represented by
an $\Bbb R$\snug-Cartier divisor $D$.  Let $U$ be the largest open subset
of $X$ such that $\pi^{-1}(U)\to U$ is an isomorphism
and that satisfies $\pi^{-1}(U)\cap\Supp D=\emptyset$.  If $\lambda$ is
a Weil function for $D$ on $W$ (defined using classical Weil functions
by $\Bbb R$\snug-linearity), then its push-forward to $U$ is
a Weil function for $\bold D$ on $X$.
\endit

If $k=\Bbb C$, then such a Weil function has domain $U(k)=U(\Bbb C)$;
if $k$ is a number field, then the domain of $\lambda$ is the disjoint
union $\coprod_{v\in M_K} U(\Bbb C_v)$, where $M_k$ is the set of places of $k$
and $\Bbb C_v$ is the completion of the algebraic closure of the local field
$k_v$ for all $v\in M_k$.

One can use such Weil functions to define proximity and counting functions
for holomorphic curves $f\:\Bbb C\to X$ whose image meets $U$ if $k=\Bbb C$,
or points in $U(k)$ or $U(\widebar k)$ if $k$ is a number field.
This is done in the obvious way.

Finally, let $Y$ be a proper closed subscheme of $X$ and let $\bold Y$ be the
corresponding b-divisor on $X$, represented by the exceptional divisor $E$
on the blowing-up $\pi\:W\to X$ of $X$ along $Y$ (as in Definition \05.7).
Then we can use this $W$ as the model in Definition \08.1, and can
use $U=X\setminus Y$.  Therefore the resulting Weil function coincides with
the Weil function as defined by \citet{sil87} or \citet{y04}
(up to an $M$\snug-bounded function, as usual).

The following is the main theorem of this section.  It generalizes
\citep{rv1, Thms.~1.4 and~1.5}, which is the special case
in which $\bold D$ is replaced by an effective (integral) Cartier divisor
on a complete variety.

\thm{\08.2}  Let $X$ be a complete variety over a field $k$, let $\Cal L$ be a
line sheaf on $X$ with $h^0(X,\Cal L^N)>1$ for some $N>0$, and let $\bold D$
be an effective $\Bbb R$\snug-Cartier b-divisor on $X$.
\roster
\myitem a. {\bc (Arithmetic part)}  Assume that $k$ is a number field,
and let $S$ be a finite set of places of $k$.  Then, for all $\epsilon>0$
and all $C\in\Bbb R$, there is a proper Zariski-closed subset $Z$ of $X$
such that the inequality
$$m_S(\bold D,x)
  \le (\Nevbir(\Cal L,\bold D)+\epsilon) h_{\Cal L}(x) + C\tag\08.2.1$$
holds for all points $x\in X(k)\setminus Z$.
\myitem b. {\bc (Analytic part)}  Assume that $k=\Bbb C$.  Then,
for all $\epsilon>0$, there is a proper Zariski-closed subset $Z$ of $X$
such that the inequality
$$m_f(\bold D,r)
  \le_\exc (\Nevbir(\Cal L,\bold D)+\epsilon) T_{f, \Cal L}(r)\tag\08.2.2$$
holds for all holomorphic mappings $f\:\Bbb C\to X$ whose image is
not contained in $Z$.
\endroster
\endit

\demo{Proof}  This is proved by reducing to the special case
\citep{rv1, Thms.~1.4 and~1.5}.

First, we reduce to the case in which $\bold D$ is replaced by an
$\Bbb R$\snug-Cartier divisor.

Let $\pi\:W\to X$ be a normal model of $X$ such that $\bold D$ is
represented by an $\Bbb R$\snug-Cartier divisor $D$ on $W$.
Since $H^0(X,\Cal L^N)\to H^0(W,\pi^{*}\Cal L^N)$ is injective,
every triple $(N,V,\mu)$ appearing in the infimum (\07.1.1) for the computation
of $\Nevbir(\Cal L,\bold D)$ also appears in the infimum for computing
$\Nevbir(\pi^{*}\Cal L,D)$.  Therefore
$$\Nevbir(\pi^{*}\Cal L,D) \le \Nevbir(\Cal L,\bold D)\;.$$
Also, as is the case for Weil functions (see Definition \08.1), we have
$$m_S(\bold D,x) = m_S(D,\pi^{-1}(x)) \qquad\text{and}\qquad
  m_f(\bold D,r) = m_{\tilde f}(D,r)$$
in the arithmetic and analytic cases, respectively, where in the arithmetic
case we consider only $x\in U$ with $U$ as in Definition \08.1,
and in the analytic case $\tilde f\:\Bbb C\to W$ is the lifting of a
holomorphic map $f\:\Bbb C\to X$ whose image meets $U$.  In both cases
we also assume that the Weil functions used for computing
the two proximity functions are related as in Definition \08.1.

For the remainder of the proof, we show only the arithmetic case;
the analytic case is similar.

By the special case for $\Bbb R$\snug-Cartier divisors, applied to $D$ on $W$,
we have
$$\split m_S(\bold D,x) &= m_S(D,\pi^{-1}(x)) \\
  &\le (\Nevbir(\pi^{*}\Cal L,D)+\epsilon) h_{\pi^{*}\Cal L}(\pi^{-1}(x)) + C \\
  &\le (\Nevbir(\Cal L,\bold D)+\epsilon) h_{\Cal L}(x) + C\;.\endsplit$$
This completes the reduction to the case of $\Bbb R$\snug-Cartier divisors.

Next, we reduce to the case in which $D$ is a $\Bbb Q$\snug-Cartier divisor.
This is done by choosing $\epsilon'\in(0,\epsilon)$ and increasing the
coefficients of $D$ by a small amount to obtain a $\Bbb Q$\snug-Cartier divisor
$D'$ such that $D\le D'\le (1+\eta)D$ with small $\eta>0$.  We then
have $m_S(D,x)\le m_S(D',x)+O(1)$.  Also
$$\split \Nevbir(\Cal L,D') + \epsilon'
  &\le \Nevbir(\Cal L,(1+\eta)D) + \epsilon' \\
  &= (1+\eta)\Nevbir(\Cal L,D) + \epsilon' \\
  &\le \Nevbir(\Cal L,D) + \epsilon\;.\endsplit$$
Here the first step is true because increasing the divisor leaves
fewer triples $(N,V,\mu)$ that satisfy (\07.1.2), which may increase the
value of the infimum; the second step holds by \citep{rv1, Remark~1.8};
and the third step is true with sufficiently small choices of $\epsilon'$
and $\eta$.

The above two inequalities then give the reduction to the case
of $\Bbb Q$\snug-Cartier divisors (apply the latter case to to $D'$).
This achieves the reduction to $\Bbb Q$\snug-Cartier divisors.

Finally, by \citep{rv1, Remark~1.8}, we can cancel the denominators
and reduce to the case of integral Cartier divisors.  This case has
been proved already \citep{rv1, Thms.~1.4 and~1.5}.\qed
\enddemo

\demo{Proof of Theorems \00.11 and \00.9}  Theorem \00.11 is immediate from
Theorem \08.2 and (\07.4.2), with $\bold D=\sum\beta(\Cal L,Y_i)\bold Y_i$.
Combining Theorem \00.11 with Theorem \00.10 (Proposition \03.3) then gives
Theorem \00.9.\qed
\enddemo

\beginsection{\09}{An Example: Linear Subspaces of $\Bbb P^n_k$}

This section gives an example involving linear subspaces of $\Bbb P^n_k$.

Let $Y_1,\dots,Y_q$ be linear subvarieties of $\Bbb P^n_k$
that intersect properly.  In this case, Definition \03.1 reduces to the
condition that they intersect properly in the sense of intersection theory;
i.e.,
$$\codim\bigcap_{i\in I}Y_i = \sum_{i\in I} \codim Y_i\tag\09.1$$
for all nonempty $I\subseteq\{1,\dots,q\}$ such that
$\bigcap_{i\in I}Y_i\ne\emptyset$.

We now compute $\beta(\Cal O(1),Y_i)$ for these subschemes.

\prop{\09.2}  Let $k$ be a field, let $X=\Bbb P^n_k$ with $n>0$, and
let $Y$ be a integral linear subscheme of $X$ of codimension $r>0$.
Then
$$\beta(\Cal O(1),Y) = \frac r{n+1}\;.\tag\09.2.1$$
\endit

\demo{Proof}  Let $x_0,\dots,x_n$ be homogeneous coordinates on $X$.
We may assume that $Y$ is the subscheme $x_1=\dots=x_r=0$.  Let $\Cal I$
be the ideal sheaf corresponding to $Y$.

We will compute $\beta(\Cal O(1),Y)$ explicitly.

First, for all $N\in\Bbb N$, $H^0(X,\Cal O(N))$ has a basis over $k$
consisting of all homogeneous monomials in $x_0,\dots,x_n$ of degree $N$.
The number of such monomials is
$$h^0(X,\Cal O(N)) = \binom{N+n}{n} = \frac{N^n}{n!} + O(N^{n-1})\tag\09.2.2$$
as $N\to\infty$, where the constant in $O(N^{n-1})$ depends only on $n$.

For all $m\in\Bbb N$ the subspace $H^0(X,\Cal O(N)\otimes\Cal I^m)$
of $H^0(X,\Cal O(N))$ is the subspace generated by
$$\{x_0^{j_0}\dotsm x_n^{j_n}
  : \text{$j_0+\dots+j_n=N$ and $j_1+\dots+j_r\ge m$}\}\;.$$
Therefore, for all $0\le m\le N$,
$$\split & h^0(X,\Cal O(N)) - h^0(X,\Cal O(N)\otimes\Cal I^m) \\
  &\qquad= \left|\left\{(j_0,\dots,j_n)\in\Bbb N^{n+1}
    : \text{$\sum j_i=N$ and $j_1+\dots+j_r<m$}\right\}\right| \\
  &\qquad= \sum_{\ell=0}^{m-1} \left|\left\{(j_0,\dots,j_n)
    : \text{$j_0+j_{r+1}+\dots+j_n=N-\ell$
      and $j_1+\dots+j_r=\ell$}\right\}\right| \\
  &\qquad= \sum_{\ell=0}^{m-1} \binom{N-\ell+n-r}{n-r}\binom{\ell+r-1}{r-1}\;.
  \endsplit\tag\09.2.3$$

For future reference, we see from the above that
$$\split & \sum_{\ell=0}^N \binom{N-\ell+n-r}{n-r}\binom{\ell+r-1}{r-1} \\
  &\qquad= \left|\left\{(j_0,\dots,j_n)
    : \text{$j_0+\dots+j_n=N$ and $j_1+\dots+j_r\le N$}
      \right\}\right| \\
  &\qquad= \binom{N+n}{n}\;.\endsplit\tag\09.2.4$$

\lemma{\09.2.5}  Let $0<r\le n$ be integers.  For all $N\in\Bbb N$ let
$$f_{n,r}(N)
  = \sum_{m=1}^N \sum_{\ell=0}^{m-1}
    \binom{N-\ell+n-r}{n-r}\binom{\ell+r-1}{r-1}\;.\tag\09.2.5.1$$
Then
$$f_{n,r}(N) = (n-r+1)\frac{N^{n+1}}{(n+1)!} + O(N^n)
  \qquad\text{as $N\to\infty$}\;,$$
where the constant in $O(N^n)$ depends only on $n$ and $r$.
\endit

\demo{Proof}  Fix $r$.  We will use induction on $n\ge r$.

For the base case $n=r$, we have
$$f_{n,n}(N) = \sum_{m=1}^N \sum_{\ell=0}^{m-1} \binom{\ell+n-1}{n-1}\;.$$
Since $\binom{\ell+n-1}{n-1}$ is a polynomial in $\ell$ with leading
term $\ell^{n-1}/(n-1)!$,
it follows that $\sum_{\ell=0}^{m-1} \binom{\ell+n-1}{n-1}$
is a polynomial in $m$ with leading term $m^n/n!$; hence $f_{n,n}(N)$ is a
polynomial in $N$ with leading term $N^{n+1}/(n+1)!$.  (In each case the
polynomial in question depends only on $n$.)

For the inductive step, assume that $n>r$ and that the lemma is true when
$n$ is replaced by $n-1$.

Since the double sum in (\09.2.5.1) is over all $m$ and $\ell$ with
$0\le\ell<m\le N$, and since the summand does not depend on $m$,
$$\split f_{n,r}(N)
  &= \sum_{\ell=0}^{N-1}
    (N-\ell) \binom{N-\ell+n-r}{n-r}\binom{\ell+r-1}{r-1} \\
  &= \sum_{\ell=0}^N
    (N-\ell) \binom{N-\ell+n-r}{n-r}\binom{\ell+r-1}{r-1}\;.\endsplit$$
By the identity $a_Nb_N-a_{N-1}b_{N-1}=a_N(b_N-b_{N-1})+(a_N-a_{N-1})b_{N-1}$,
Pascal's Rule, (\09.2.4), the inductive hypothesis, and (\09.2.2),
$$\split f_{n,r}(N) - f_{n,r}(N-1)
  &= \sum_{\ell=0}^{N-1} (N-\ell)\binom{N-\ell+n-r-1}{n-r-1}
    \binom{\ell+r-1}{r-1} \\
  &\qquad+ \sum_{\ell=0}^{N-1} \binom{N-\ell+n-r-1}{n-r}\binom{\ell+r-1}{r-1} \\
  &= f_{n-1,r}(N) + \binom{N-1+n}{n} \\
  &= \left((n-r)\frac{N^n}{n!} + O(N^{n-1})\right)
    + \left(\frac{(N-1)^n}{n!} + O((N-1)^{n-1})\right) \\
  &= (n-r+1)\frac{N^n}{n!} + O(N^{n-1})\endsplit$$
as $N\to\infty$, where again the implicit constants depend only on $n$ and $r$.

Since $f_{n,r}(0)=0$, the lemma then follows by induction.\qed
\enddemo

Applying (\09.2.3), the lemma, and (\09.2.2), we then have
$$\split \beta(\Cal O(1),Y)
  &= \liminf_{N\to\infty}
    \frac{\sum_{m=1}^N h^0(X,\Cal O(N)\otimes\Cal I^m)}{Nh^0(X,\Cal O(N))} \\
  &= \liminf_{N\to\infty}
    \frac{Nh^0(X,\Cal O(N)) - f_{n,r}(N)}{Nh^0(X,\Cal O(N))} \\
  &= 1 - \limsup_{N\to\infty}
    \frac{(n-r+1)N^{n+1}/(n+1)!+O(N^n)}{N\cdot N^n/n!+O(N^n)} \\
  &= 1 - \frac{n-r+1}{n+1} \\
  &= \frac r{n+1}\;,\endsplit$$
as was to be shown.\qed
\enddemo

As a corollary of Theorem \00.9, we then obtain:

\thm{\09.3}  Let $k$ be a number field, let $S$ be a finite set of places
of $k$, let $X=\Bbb P^n_k$, and let $Y_1,\dots,Y_q$ be
linear subvarieties of $X$ in general position (according to (\09.1)).
Then, for all $\epsilon>0$ and all $C\in\Bbb R$, there is a
proper Zariski-closed subset $Z$ of $X$ such that the inequality
$$\sum_{i=1}^q (\codim Y_i)m_S(Y_i,x) \le (n+1+\epsilon)h_k(x) + C$$
holds for all $x\in X(k)$ outside of $Z$.
\endit

This is a consequence of \citep{voj89, (3.9)}, in which $\Cal D$ is a
finite collection of hyperplanes containing, for each $i$, a subset whose
intersection is $Y_i$.

It is also a special case of the Main Theorem of \citep{hl17}.

\beginsection{\010}{Proof that (\00.1.1), (\00.4.1), (\00.5.1), and (\05.5.1)
  Exist as Limits}

This section gives a proof that the limits infima in the definitions of
$\beta(\Cal L,D)$ (Definition \00.1), $\beta(\Cal L,Y)$ (Definition \00.4),
$\beta_{\Cal L,Y}$ (Remark \00.5), and
$\beta(\Cal L,\bold D)$ (Definition \05.5)
can be replaced by limits.

It has already been noted that the lim inf in the definition of
$\beta(\Cal L,D)$ (Definition \00.1) is a limit when $D$ is a
nonzero effective Cartier divisor (see the discussion following Def.~1.9
in \citep{rv1}).
We extend this result to allow $D$ to be a nonzero effective Cartier b-divisor.
Since a detailed proof has not appeared before, we include here such a proof
of both results.  It will then be immediate from Corollary \05.9b that the
same is true for $\beta(\Cal L,Y)$ and $\beta_{\Cal L,Y}$.

Recall that in all cases, $\Cal L$ is assumed to be big.

This argument is based on an idea of Julie Wang to compare the limit
to a Riemann sum.

We start with the proof that the limit in (\00.1.1) converges.

\thm{\010.1}  Let $X$ be a complete variety over a field $F$
of characteristic zero, let $\Cal L$ be a big line sheaf on $X$,
and let $D$ be a nonzero effective Cartier divisor on $X$.  Then the limit
$$\lim_{N\to\infty}
    \frac{\sum_{m=1}^\infty h^0(X,\Cal L^N(-mD))}{Nh^0(X,\Cal L^N)}
  \tag\010.1.1$$
converges.  In particular, the $\liminf$ in Definition \00.1 can be
replaced by a limit.
\endit

\demo{Proof}  We start by reducing to the projective case.  Let $d=\dim X$.

By Chow's lemma and resolution of singularities there is a model $\pi\:W\to X$,
with $W$ projective and nonsingular.  Then
$$\sum_{m=1}^\infty h^0(X,\Cal L^N(-mD))
  = \sum_{m=1}^\infty h^0(W,\pi^{*}(\Cal L^N(-mD))) + O(N^d)$$
by Lemma \05.2d, and
$$h^0(X,\Cal L^N) = h^0(W,\pi^{*}\Cal L^N) + O(N^{d-1})$$
by Lemma \05.2c.  Therefore
$$\lim_{N\to\infty}
    \frac{\sum_{m=1}^\infty h^0(X,\Cal L^N(-mD))}{Nh^0(X,\Cal L^N)}
  = \lim_{N\to\infty}
    \frac{\sum_{m=1}^\infty h^0(W,\pi^{*}(\Cal L^N(-mD)))}
      {Nh^0(W,\pi^{*}\Cal L^N)}\;,$$
in the sense that if one limit converges, then both do, and they are equal.

So assume now that $X$ is projective and nonsingular.

For all line sheaves $\Cal L$ on $X$, all effective Cartier divisors
$D$ on $X$, and all $x\in\Bbb R_{\ge0}$, we let
$$H^0(X,\Cal L(-xD)) = \{s\in H^0(X,\Cal L):
  \text{the $\Bbb R$\snug-divisor $(s)-xD$ is effective}\}$$
and (as usual)
$$h^0(X,\Cal L(-xD)) = \dim_F H^0(X,\Cal L(-xD))\;.$$
These coincide with the usual definitions whenever $xD$ is an integral divisor.

Define $f\:\Bbb R_{\ge0}\to\Bbb R$ by
$$f(x) = \lim_{N\to\infty} \frac{h^0(X,\Cal L^N(-NxD))}{N^d}\;,$$
where the limit is over $N\in\Bbb Z_{>0}$.

Recall from \citep{laz04, II, Def.~11.4.2 and Example 11.4.7} that
$$\lvol(\Cal L) = \lim_{N\to\infty} \frac{h^0(X,\Cal L^N)}{N^d/d!}\;.
  \tag\010.1.2$$
Then $f(x)=\lvol(\Cal L(-xD))/d!$ whenever $xD$ is an integral Cartier divisor.

Since $D$ is effective, $f$ is a nonincreasing function.

We also have $f(x)=0$ for all sufficiently large $x$.  Indeed,
given an ample divisor $A$ on $X$, this is true for all
$x>(\Cal L\idot A^{d-1})/(D\idot A^{d-1})$.  Fix some $R\in\Bbb R_{\ge0}$
such that $H^0(X,\Cal L^N(-NRD))=0$ for all $N>0$ (and therefore $f(R)=0$).

Let
$$I = \int_0^\infty f(x)\,dx = \int_0^R f(x)\,dx\;.$$
It will then suffice to prove that
$$\lim_{N\to\infty} \frac1{N^{d+1}}\sum_{m\ge1} h^0(X,\Cal L^N(-mD)) = I\;,
  \tag\010.1.3$$
since by (\010.1.2) this would imply
$$\lim_{N\to\infty} \frac{\sum_{m\ge1} h^0(X,\Cal L^N(-mD))}{Nh^0(X,\Cal L^N)}
  = \frac {d!I}{\lvol(\Cal L)}\;.$$

\lemma{\010.1.4}  Let $\Cal M$ be a line sheaf on $X$.  If the limit
$$\lim_{N\to\infty} \frac1{N^{d+1}}\sum_{m\ge1}h^0(X,\Cal L^N(-mD))$$
exists, then so does the limit
$$\lim_{N\to\infty} \frac1{N^{d+1}}\sum_{m\ge1}h^0(X,\Cal M\otimes\Cal L^N(-mD))
  \;,$$
and they are equal.
\endit

\demo{Proof}  First, let $\Cal M_1$ and $\Cal M_2$ be line sheaves on $X$.
Since $X$ is projective, we have
$h^0(X,\Cal M_1\otimes\Cal M_2^{-1}\otimes\Cal L^p)\ne0$ for some $p\in\Bbb N$
by a consequence of Kodaira's lemma \citep{laz04, 2.2.7}.  Therefore
$$h^0(X,\Cal M_2\otimes\Cal L^N(-mD))
  \le h^0(X,\Cal M_1\otimes\Cal L^{N+p}(-mD))$$
for all $m,N\in\Bbb Z_{>0}$.  Therefore we have
$$\split & \liminf_{N\to\infty}
    \frac1{N^{d+1}}\sum_{m\ge1} h^0(X,\Cal M_2\otimes\Cal L^N(-mD)) \\
  &\qquad\le \liminf_{N\to\infty}
    \frac1{(N-p)^{d+1}}\sum_{m\ge1} h^0(X,\Cal M_1\otimes\Cal L^N(-mD)) \\
  &\qquad= \liminf_{N\to\infty}
    \frac1{N^{d+1}}\sum_{m\ge1} h^0(X,\Cal M_1\otimes\Cal L^N(-mD))\;,
  \endsplit$$
and likewise for $\limsup$.

This gives
$$\split \liminf_{N\to\infty} \frac1{N^{d+1}} \sum_{m\ge1} h^0(X,\Cal L^N(-mD))
  &\le \liminf_{N\to\infty}
    \frac1{N^{d+1}} \sum_{m\ge1} h^0(X,\Cal M\otimes\Cal L^N(-mD)) \\
  &\le \limsup_{N\to\infty}
    \frac1{N^{d+1}} \sum_{m\ge1} h^0(X,\Cal M\otimes\Cal L^N(-mD)) \\
  &\le \limsup_{N\to\infty} \frac1{N^{d+1}} \sum_{m\ge1} h^0(X,\Cal L^N(-mD))\;,
  \endsplit$$
and this implies the lemma.\qed
\enddemo

Now let $k\in\Bbb Z_{>0}$.  We show that if
$$\lim_{N\to\infty}
    \frac1{N^{d+1}}\sum_{m\ge1} h^0(X,\Cal M\otimes\Cal L^{Nk}(-mD))
  = k^{d+1}I\tag\010.1.5$$
with $\Cal M=\Cal O_X$, then (\010.1.3) is true.  Indeed, if (\010.1.5)
is true with $\Cal M=\Cal O_X$, then by Lemma \010.1.4
it is true with $\Cal M=\Cal L^j$ with $j=0,1,\dots,k-1$.
Therefore the limit in (\010.1.3) exists for $N$ in each congruence class
modulo $k$, and these limits are all equal.

Thus, if (\010.1.3) is true with $\Cal L$ replaced by $\Cal L^k$
for some $k>0$, then it is true with the original $\Cal L$.
In particular, choosing $k$ such that $H^0(X,\Cal L^k)\ne0$, we may
assume that $H^0(X,\Cal L)\ne0$.

We then have
$$h^0(X,\Cal L^N(-mD)) \le h^0(X,\Cal L^{N'}(-mD))\tag\010.1.6$$
for all $0\le N\le N'$ and all $m\in\Bbb N$.

We now begin the main argument of the proof.

Given $\epsilon>0$, pick $\epsilon_1>0$ and $k,l_0\in\Bbb Z_{>0}$ such that
$$\left(1+\frac1{l_0}\right)^{d+1}\left(I+\frac{f(0)}{k}+\epsilon_1\right)
  \le I+\epsilon\tag\010.1.7$$
and
$$\left(1-\frac1{l_0}\right)^{d+1}\left(I-\frac{f(0)}{k}-\epsilon_1\right)
  \ge I-\epsilon\;.$$

We claim that if $k$ and $l_0$ are chosen sufficiently large, then we also have
$$\frac1{(lk)^d}\sum_{m=0}^\infty h^0(X,\Cal L^{lk}(-mlD))
  \le \sum_{m=0}^\infty f\fracwithdelims()mk + \epsilon_1k\tag\010.1.8$$
and
$$\frac1{(lk)^d}\sum_{m=1}^\infty h^0(X,\Cal L^{lk}(-mlD))
  \ge \sum_{m=1}^\infty f\fracwithdelims()mk - \epsilon_1k\tag\010.1.9$$
for all $l\ge l_0$.

We will show this result only for (\010.1.8).  The argument for (\010.1.9)
is similar and is left to the reader.

We may assume that $R\in\Bbb Z$.

Choose $\epsilon_2>0$, $\epsilon_3>0$, and $\epsilon_4>0$
such that $\epsilon_2+\epsilon_3+\epsilon_4\le\epsilon_1$.

Choose $x_0,\dots,x_t\in\Bbb R$ such that $0=x_0<x_1<\dots<x_t=R$ and
$$\sum_{i=1}^t (x_i-x_{i-1})(f(x_{i-1})-f(x_i)) \le \epsilon_2\;.$$
Define a function $g\:\Bbb R_{\ge0}\to\Bbb R_{\ge0}$ by
$$g(x) = f(x_{i-1})-f(x)
  \qquad\text{for all $x\in\sq(x_{i-1},x_i)$ and all $i$}$$
and by $g(x)=0$ for all $x\ge R$.  Then $g(x)\le f(x_{i-1})-f(x_i)$
for all $x\in\sq(x_{i-1},x_i)$ and all $i$; hence
$$\int_0^\infty g(x)\,dx = \int_0^R g(x)\,dx \le \epsilon_2\;.$$
By the theory of Riemann integration, since $g$ is piecewise nondecreasing,
we have
$$\frac1k \sum_{m=0}^\infty g\fracwithdelims()mk \le \epsilon_2+\epsilon_3$$
for all sufficiently large $k$.  Fix such a $k$.  Then there is an
integer $N_0$, depending on $k$, such that
$$\left|\frac{h^0(X,\Cal L^N(-Nx_{i-1}D))}{N^d} - f(x_{i-1})\right|
  \le \frac{\epsilon_4}{kR}$$
for all $1\le i\le t$ and all $N\ge N_0$.  Therefore
$$\frac{h^0(X,\Cal L^N(-NxD))}{N^d} \le \frac{h^0(X,\Cal L^N(-Nx_{i-1}D))}{N^d}
  \le f(x_{i-1}) + \frac{\epsilon_4}{kR}
  = f(x) + g(x) + \frac{\epsilon_4}{kR}$$
for all $i$, all $x\in[x_{i-1},x_i]$, and all $N\ge N_0$.
Let $l_0=\lceil N_0/k\rceil$.  Then, for all $l\ge l_0$,
$$\split & \frac1k \sum_{m=0}^\infty \frac{h^0(X,\Cal L^{lk}(-mlD))}{(lk)^d}
    - \frac1k \sum_{m=0}^\infty f\fracwithdelims()mk \\
  &\qquad= \frac1k \sum_{m=0}^{kR-1} \frac{h^0(X,\Cal L^{lk}(-mlD))}{(lk)^d}
    - \frac1k \sum_{m=0}^{kR-1} f\fracwithdelims()mk \\
  &\qquad\le \frac1k \sum_{m=0}^{kR-1} g\fracwithdelims()mk + \epsilon_4 \\
  &\qquad\le \epsilon_2+\epsilon_3+\epsilon_4 \\
  &\qquad\le \epsilon_1\;.
  \endsplit$$
This concludes the proof of the claim.

By elementary facts about Riemann sums for monotone functions, we have
$$\frac1k\sum_{m=1}^{kR} f\fracwithdelims()mk \le I
  \le \frac1k\sum_{m=0}^{kR-1} f\fracwithdelims()mk\;.$$
Since $f(R)=0$ and since the two sums differ by $f(0)/k$, we have
$$\frac1k\sum_{m=0}^{kR} f\fracwithdelims()mk \le I + \frac{f(0)}{k}
  \qquad\text{and}\qquad
  \frac1k\sum_{m=1}^{kR} f\fracwithdelims()mk \ge I - \frac{f(0)}{k}\;.
  \tag\010.1.10$$

Now let any $N\ge l_0k$ be given.  Let $l=\fracwithdelims\lceil\rceil Nk$.
Then $lk\ge N$ and $l\ge l_0$; hence
$$\frac{lk}{N} < \frac{N+k}{N} \le 1 + \frac1{l_0}\;.\tag\010.1.11$$
By (\010.1.6), effectivity of $D$, (\010.1.11), (\010.1.8), (\010.1.10),
and (\010.1.7),
$$\split \frac1{N^{d+1}} \sum_{m\ge1} h^0(X,\Cal L^N(-mD))
  &\le \frac1{N^{d+1}} \sum_{m\ge1} h^0(X,\Cal L^{lk}(-mD)) \\
  &\le \frac1{N^{d+1}} \sum_{m\ge0} h^0\left(X,
    \Cal L^{lk}\left(-\fracwithdelims\lfloor\rfloor ml lD\right)\right) \\
  &= \frac 1{N^{d+1}} \sum_{m'=0}^\infty lh^0(X,\Cal L^{lk}(-m'lD)) \\
  &< \left(1+\frac1{l_0}\right)^{d+1}
    \frac1{l^dk^{d+1}}\sum_{m=0}^\infty h^0(X,\Cal L^{lk}(-mlD)) \\
  &\le \left(1+\frac1{l_0}\right)^{d+1}
    \left(\frac1k \sum_{m=0}^\infty f\fracwithdelims()mk + \epsilon_1\right) \\
  &\le \left(1+\frac1{l_0}\right)^{d+1}
    \left(I + \frac{f(0)}{k} + \epsilon_1\right) \\
  &\le I + \epsilon\;.\endsplit$$
A similar argument gives
$$\frac1{N^{d+1}} \sum_{m\ge1} h^0(X,\Cal L^N(-mD)) > I - \epsilon\;,$$
and this implies (\010.1.3), concluding the proof of the theorem.\qed

\cor{\010.2}  Let $X$ be a complete variety over a field of characteristic $0$,
let $\Cal L$ be a big line sheaf on $X$, and let $\bold D$
be a nonzero effective Cartier b-divisor on $X$.
Then, recalling Definition \05.3, the limit
$$\lim_{N\to\infty}
    \frac{\sum_{m=1}^\infty h^0_\bir(X,\Cal L^N(-m\bold D))}
      {Nh^0_\bir(X,\Cal L^N)}\tag\010.2.1$$
exists.  Thus, the $\liminf$ in (\05.5.1) is actually a limit.
\endit

\demo{Proof}  Let $\pi\:W\to X$ be a normal model of $X$ on which $\bold D$
is represented by an effective Cartier divisor $D$.
By Lemma \05.2a, we then have
$$h^0_\bir(X,\Cal L^N(-m\bold D)) = h^0(W,\pi^{*}\Cal L^N(-mD))$$
for all $m,N\in\Bbb N$ (notably including $m=0$).  Thus
$$\frac{\sum_{m=1}^\infty h^0_\bir(X,\Cal L^N(-m\bold D))}
    {Nh^0_\bir(X,\Cal L^N)}
  = \frac{\sum_{m=1}^\infty h^0(W,\pi^{*}\Cal L^N(-mD))}
    {Nh^0(W,\pi^{*}\Cal L^N)}\qquad\text{for all $N\in\Bbb Z_{>0}$}\;,$$
and the corollary then follows from Theorem \010.1.\qed
\enddemo

\cor{\010.3}  If $\operatorname{char} k=0$, then the limits infima
in (\00.4.1) and (\00.5.1) converge as limits.
\endit

\demo{Proof}  This is immediate from Corollaries \010.2 and \05.9b.\qed
\enddemo

\comment
Here's the ``similar argument'':

Let any $N\ge l_0k$ be given.  Let $l=\fracwithdelims\lfloor\rfloor Nk$.
Then $lk\le N$ and $l\ge l_0$; hence
$$\frac{lk}{N} > \frac{N-k}{N} \ge 1 - \frac1{l_0}\;.\tag\010.1.12$$
By (\010.1.6), effectivity of $D$, (\010.1.12), (\010.1.8), (\010.1.10),
and (\010.1.7),
$$\split \frac1{N^{d+1}} \sum_{m\ge1} h^0(X,\Cal L^N(-mD))
  &\ge \frac1{N^{d+1}} \sum_{m\ge1} h^0(X,\Cal L^{lk}(-mD)) \\
  &\ge \frac1{N^{d+1}} \sum_{m\ge1} h^0\left(X,
    \Cal L^{lk}\left(-\fracwithdelims\lceil\rceil ml lD\right)\right) \\
  &= \frac l{N^{d+1}} \sum_{m'=1}^\infty h^0(X,\Cal L^{lk}(-m'lD)) \\
  &> \left(1-\frac1{l_0}\right)^{d+1}
    \frac1{l^dk^{d+1}}\sum_{m=1}^\infty h^0(X,\Cal L^{lk}(-mlD)) \\
  &> \left(1-\frac1{l_0}\right)^{d+1}
    \left(\frac1k \sum_{m=1}^\infty f\fracwithdelims()mk - \epsilon_1\right) \\
  &\ge \left(1-\frac1{l_0}\right)^{d+1}
    \left(I - \frac{f(0)}{k} - \epsilon_1\right) \\
  &\ge I - \epsilon\;.\endsplit$$

\endcomment

\enddemo

\Refs

\newdimen\bibindent
\setbox0\hbox{--1234pre--}
\bibindent=\wd0

\newdimen\bibitemindent
\setbox0\hbox{--}
\bibitemindent=\wd0

\auth{Autissier, Pascal}
\refer{a} Sur la non-densit\'e des points entiers.
  {\it Duke Math. J.} 158 (2011) pp.~13--27.

\auth{Hartshorne, Robin}
\refer{ha77} {\it Algebraic Geometry,} Graduate Texts in Mathematics, No. 52,
  Springer-Verlag, New York-Heidelberg, 1977.

\auth{Heier, Gordon and Levin, Aaron}
\refer{hl17} A generalized Schmidt subspace theorem for closed subschemes.
  Preprint arXiv:1712.02456 (2017).

\auth{Jothilingham, P., Duraivel, T., and Bose, Sibnath}
\refer{JDB} Regular sequences of ideals.
  {\it Beitr\. Algebra Geom.} 52 (2011), pp.~479--485.

\auth{Lang, Serge}
\refer{la_fdg}  {\it Fundamentals of {D}iophantine geometry,}
  Springer-Verlag, New York, 1983.

\auth{Lazarsfeld, Robert}
\refer{laz04} {\it Positivity in Algebraic Geometry I, II,}
  Ergebnisse der Mathematik und ihrer Grenzgebiete 3. Folge 48, 49,
  Springer-Verlag, Berlin Heidelberg 2004.

\auth{Ru, Min and Vojta, Paul}
\refer{rv1} A birational Nevanlinna constant and its consequences.
  {\it Amer. J. Math.} 142 (2020), pp.~957--991.

\refer{rv2} An Evertse--Ferretti Nevanlinna constant and its consequences.
  Preprint arXiv:2004.12257 (2020).

\auth{Ru, Min and Wang, Julie Tzu-Yueh}
\refer{rw17} A subspace theorem for subvarieties.
  {\it Algebra and Number Theory} 11 (2017), pp.~2323--2337.

\refer{rw20p}  The Ru--Vojta result for subvarieties.
  Preprint (2020).

\auth{Silverman, Joseph H.}
\refer{sil87} Arithmetic distance functions and height functions in diophantine
  geometry.
  {\it Math. Ann.} 279 (1987) pp.~193--216.

\auth{The Stacks project authors}
\refer{stks} The Stacks project,
  {\tt https://stacks.math.columbia.edu} (2020).

\auth{Vojta, Paul}
\refer{voj89} A refinement of Schmidt's Subspace Theorem.
  {\it Amer. J. Math.} 111 (1989), pp.~489--518.

\refer{voj11} Diophantine approximation and Nevanlinna theory.
  In: {\it Arithmetic geometry,} Lecture Notes in Mathematics 2009,
  Springer, Berlin, 2011, pp.~111--224.

\auth{Yamanoi, Katsutoshi}
\refer{y04} Algebro-geometric version of Nevanlinna's lemma on
  logarithmic derivative and applications.
  {\it Nagoya Math. J.} 173 (2004) pp.~23--63.

\endRefs

\enddocument